\documentclass[11pt,letterpaper,reqno]{amsart} 
\usepackage{graphicx} 
\usepackage{amsfonts, amsmath, amssymb, amsthm, multicol, enumitem, mathtools, bm}
\usepackage[left=1in, right=1in, bottom=1in, top=1in, bottom=1in]{geometry}

\usepackage[dvipsnames]{xcolor}

\usepackage{amsrefs}
\usepackage{tikz}
\usetikzlibrary{positioning}
\usetikzlibrary{arrows.meta}

\usepackage{etoolbox}
\usepackage{hycolor} 
\usepackage[colorlinks=true]{hyperref} 
\hypersetup{
  colorlinks=true,
  linkcolor=black,
  citecolor=blue,
  linkbordercolor={0 0 1}
} 
\usepackage{cleveref} 

\title{K3 surfaces and BHK Mirror symmetry: A counterexample}
\author{Alexis Barrett}
\author{Archer Clayton}
\author{Ben DeVries}
\author{Nathan Priddis}

\date \today

\newtheorem{theorem}{Theorem}[section]

\newtheorem{conjecture}{Conjecture}
\newtheorem{lemma}[theorem]{Lemma}

\theoremstyle{definition}
\newtheorem{definition}[theorem]{Definition}

\DeclareMathOperator{\SL}{SL}
\DeclareMathOperator{\GL}{GL}
\DeclareMathOperator{\Aut}{Aut}
\DeclareMathOperator{\Hom}{Hom}
\DeclareMathOperator{\Pic}{Pic}
\DeclareMathOperator{\image}{Im}

\DeclareMathOperator{\rk}{rk}

\newcommand{\PP}{\mathbb{P}}
\newcommand{\CP}{\mathbb{CP}}
\newcommand{\ZZ}{\mathbb{Z}}
\newcommand{\CC}{\mathbb{C}}
\newcommand{\QQ}{\mathbb{Q}}

\newcommand{\Gmax}{G^{\text{max}}}
\newcommand{\Gdiag}{G^{\text{diag}}}
\newcommand{\aplus}{a_+}
\newcommand{\aminus}{a_-}
\newcommand{\bplus}{b_+}
\newcommand{\bminus}{b_-}

\begin{document}

\begin{abstract}
For K3 surfaces of BHK type, there are two versions of mirror symmetry that might be applied, namely BHK mirror symmetry and LPK3 mirror symmetry of Dolgachev, which involved choosing a polarizing lattice. It has been conjectured that for a pair of BHK mirror K3 surfaces, one can choose a polarizing lattice as the invariant lattice of some nonsymplectic diagonal symmetry. 
We provide a counterexample, showing that LPK3 mirror symmetry cannot be realized by polarizing surfaces by the invariant lattice of a diagonal symmetry.  However, in our example the polarizing lattice can be realized as the invariant lattice of a permutation symmetry. 
\end{abstract}

\maketitle

\tableofcontents

\section{Introduction}

Originating from string theory, mirror symmetry is a remarkable duality between families of Calabi–Yau manifolds of the same dimension. Physicists observed that seemingly different Calabi–Yau manifolds can define equivalent physical theories. Since its emergence roughly three decades ago, mirror symmetry has led to a wide range of mathematical predictions and discoveries.

Mathematically, mirror symmetry asserts that with a Calabi–Yau manifold $X$ one can associate a mirror Calabi-Yau manifold $X^\vee$ whose complex geometry reflects the Kähler geometry of $X$, and vice versa. More precisely, the structures governing the Kähler deformations of $X$ are expected to correspond, in an appropriate sense, to the structures governing the complex deformations of $X^\vee$.

Mathematicians have devised several different constructions, depending on the context, that describe mirror symmetry. In this article, we consider two-dimensional Calabi-Yau manifolds called K3 surfaces. We focus on a special class of K3 surfaces that arise as quotients of hypersurfaces in weighted projective spaces. Such K3 surfaces, denoted $X_{W,G}$, are determined by a pair $(W,G)$ of a defining polynomial $W$ and a group of symmetries $G$ used to make the quotient. These K3 surfaces lie in the intersection of two mirror symmetry constructions.

The two forms of mirror symmetry that we consider are known as BHK mirror symmetry and LPK3 mirror symmetry. BHK mirror symmetry defines a mirror pair $(W^T,G^T)$ used to construct the mirror K3 surface to the surface defined by the original pair $(W,G)$. In contrast, LPK3 mirror symmetry involves polarizing a K3 surface $X$ by an even integral lattice $M$. If $X^\vee$ can be polarized by the mirror lattice $M^\vee$, then we say $X$ and $X^\vee$ are LPK3 mirror duals. We describe each of these in more detail in Section~\ref{sec:mirror-symmetry}.

The other form of mirror symmetry, LPK3 mirror symmetry, involves polarizing a K3 surface $X$ by an even integral lattice $M$. If $X^\vee$ can be polarized by the mirror lattice $M^\vee$, then we say $X$ and $X^\vee$ are LPK3 mirror duals.

The natural question is whether these two versions of mirror symmetry ``agree'' for K3 surfaces $X_{W,G}$ defined by a pair $(W,G)$. In other words, for every BHK mirror pair, we would like to know if there exists some polarization of each surface such that the two surfaces are also LPK3 mirror duals. This article addresses this question for three K3 surfaces.

\subsection{Previous results}

The question of whether BHK mirror symmetry and LPK3 mirror symmetry agree has been studied in various settings by several authors.  
About 15 years ago, in \cite{ABS}, Artebani, Boissierre, and Sarti proved that the two forms of mirror symmetry do ``agree'' for K3 surfaces defined by a polynomial of the form $x^2+f(y,z,w)$. Later results by several other authors (see \cite{BCP}, \cite{CLPS}, and \cite{CP}) extended this first result to polynomials of the form $x^n+f(y,z,w)$ as follows (see Section~\ref{sec:background} for definitions):  

\begin{theorem} \label{thm:fermat-bhk-lpk3-mirror}
Let $W$ be an invertible polynomial of the form 
\begin{equation}\label{eq:x^nform}
      W=x^n+f(y,z,w),
\end{equation}
and let $G$ be a group of diagonal symmetries satisfying $J_W\subseteq G\subseteq \SL_W$, so that $X_{W,G}$ is a K3 surface. Then $X_{W^T,G^T}$ is an LPK3 dual to $X_{W,G}$.  
\end{theorem}

Due to the form \eqref{eq:x^nform} of the polynomial $W$, there is a canonical purely nonsymplectic automorphism $\sigma_n$ in $X_{W,G}$ given by 
\[ 
\sigma_n: [x:y:z:w]\mapsto [\xi_nx:y:z:w].
\]
As the dual polynomial will also be of the same form \eqref{eq:x^nform}, there is a similar purely nonsymplectic automorphism $\sigma_n^T$ on $X_{W^T,G^T}$. Theorem~\ref{thm:fermat-bhk-lpk3-mirror} is proved by polarizing $X_{W,G}$ by the invariant lattice $S(\sigma_n)$ of $\sigma_n$ and by polarizing $X_{W^T,G^T}$ by the invariant lattice $S(\sigma_n^T)$ of $\sigma_n^T$. 

The proof is completed by showing that the invariant lattice $S(\sigma_n^T)$ is the mirror of $S(\sigma_n)$. That is, given two K3 surfaces of the form \eqref{eq:x^nform} that are BHK mirrors, one must find a specific nonsymplectic automorphism $\sigma_n$ on $X_{W,G}$, and a \emph{specific} nonsymplectic automorphism $\sigma_n^T$ on $X_{W^T,G^T}$ so that the respective invariant lattices satisfy $S(\sigma_n^T)\cong S(\sigma_n)^\vee$.

It is well-known that there are 95 weight systems that produce K3 surfaces of the form $X_{W,G}$ when $W$ is nondegenerate and quasihomogeneous with respect to one of these weight systems. However, Theorem~\ref{thm:fermat-bhk-lpk3-mirror} is only known for invertible polynomials of the form \eqref{eq:x^nform}. There are six weight systems that admit invertible polynomials, but do not admit polynomials of this form. Moreover, nearly every one of the 95 weight systems admits a polynomial not of this form. Moreover, each polynomial generates several different K3 surfaces due to differing choices of the symmetry group $G$, and thus different BHK mirror pairs. 

This naturally leads us to the following Conjecture. 
\begin{conjecture}\label{conj:main}
  Let $X_{W,G}$ be a K3 surface determined by the pair $(W,G)$, where $W$ is invertible with respect to one of the 95 weight systems, and $J_W\subseteq G\subseteq \SL_W$. Then there exists a (diagonal) nonsymplectic automorphism $\sigma$ on $X_{W,G}$ and a (diagonal) nonsymplectic automorphism $\sigma^\vee$ on $X_{W^T,G^T}$ so that $S(\sigma)^\vee\cong S(\sigma^\vee)$. In other words, the polarizing lattice for $X_{W,G}$ and its mirror lattice polarizing $X_{W^T,G^T}$ are invariant lattices for some nonsymplectic diagonal symmetries. 
\end{conjecture}

We focus on diagonal symmetries for several reasons. First, 
diagonal automorphisms are simple to work with because their action can be understood explicitly in terms of (weighted) homogeneous coordinates in weighted projective space. Second, they appear to capture much of the geometry of a K3 surface. For example, in \cite{Br} Brandhorst proves that K3 surfaces with a purely nonsymplectic automorphism of a particular order are unique up to isomorphism. It seems that all of these automorphisms can be realized as a diagonal symmetry for some K3 surface of the form $X_{W,G}$ (this is a pending result of the last author). One might also notice that each of the automorphisms used to prove Theorem~\ref{thm:fermat-bhk-lpk3-mirror} are also diagonal symmetries. 

We investigate three K3 surfaces constructed from a polynomial $W$ that is not of the form \eqref{eq:x^nform}, hoping to verify or disprove Conjecture~\ref{conj:main}. These three surfaces give us two ``pairs'' of BHK mirror surfaces---two are mirror dual, and the third is self-dual.

\subsection{Main Results}
The first main result is that Conjecture~\ref{conj:main} is not true for \textit{diagonal symmetries}. For each of the three K3 surfaces in question, we investigate the invariant lattices of every possible nonsymplectic diagonal automorphism\footnote{We wrote Python code to verify all our computations. This code can be found in the Github repository for this paper \cite{DeVries}} to verify that no invariant lattice of such an automorphism on the first surface is the mirror of an invariant lattice on the dual surface.

For a second result, We also consider certain nonsymplectic automorphisms that are not diagonal symmetries, and in every case, we find a polarizing lattice that makes each BHK mirror pair also into an LPK3 mirror pair. However, the automorphism is very different in nature from the automorphism used to prove Theorem~\ref{thm:fermat-bhk-lpk3-mirror}. 

Each K3 surface that we consider has Picard rank 20, and every K3 surface with Picard rank 20 has an infinite automorphism group (see e.g. \cite{ShI} or \cite[Chapter 15]{huybrechtsK3}). Our work leads to the following natural question. Can one determine if Conjecture~\ref{conj:main} is true if one does not insist upon using a diagonal symmetry, and if so, is there some canonical choice of nonsymplectic automorphism that provides an appropriate polarizing lattice in every case.

\subsection{Outline}
In Section~\ref{sec:background}, we will go over the requisite background for this work, including defining K3 surfaces, discussing conditions on $W$ and $G$, introducing integral lattices, and finally introducing two forms of mirror symmetry for K3 surfaces. In Section~\ref{sec:example}, we describe three K3 surfaces that are constructed using the same polynomial $W$, but using three different quotient groups. We describe the geometry, the diagonal symmetries, and discuss the LPK3 mirror symmetry for each pair of BHK mirrors. As the majority of the authorship of this article consists of students, our target audience is advanced undergraduate students and graduate students. Nevertheless we hope anyone who reads this learns an interesting result.

\subsection{Acknowledgments} The authors would like to thank Brigham Young University (BYU) for its generous support of undergraduate research. Three of the authors were undergraduates at BYU during the majority of this project. They were supported by generous donors to the BYU Math Department and BYU's College of Computational, Mathematical, and Physical Sciences. 

Use of AI on this project was very limited. ChatGPT was consulted to learn several references, and Claude was used to assist creating some of the diagrams.

\section{Background} \label{sec:background}

\subsection{K3 surfaces and automorphisms} \label{sec:k3-surfaces-and-automorphisms}

Most of the definitions and facts in this section on K3 surfaces are well-known. For a reference, see e.g. \cite{huybrechtsK3}. 

\begin{definition}
    A \emph{K3 surface} is a compact complex surface that is simply connected and has a trivial canonical bundle. 
\end{definition}

It is well-known that all K3 surfaces are diffeomorphic. Given a K3 surface $X$, the cohomology $H^2(X,\ZZ)$ is a free abelian group with rank 22, and the Euler characteristic of $X$ is 24. The nonzero Hodge numbers are $h^{0,0}(X)=h^{2,2}(X)=h^{0,2}(X)=h^{2,0}(X)=1$, and $h^{1,1}(X)=20$. This gives us the following Hodge diamond: 
\begin{equation}\label{eq:K3HodgeDiamond}
\begin{matrix}
     & 1 & \\
     1 & 20 & 1\\\
      & 1 & \\
\end{matrix}    
\end{equation}

Because the canonical bundle is trivial and $h^{2,0}(X)=1$, every K3 surface $X$ has a non-vanishing volume form $\omega$ that generates $H^{2,0}(X)$. Hence, any automorphism $\sigma\in \Aut(X)$ induces an action on $\omega$, via $\sigma^*\omega=\lambda_\sigma\omega$ for some $\lambda_\sigma\in \CC^*$. If $\lambda_\sigma=1$, then we say $\sigma$ is \emph{symplectic}. Otherwise, $\sigma$ is \emph{nonsymplectic}. If $\sigma$ has finite order $n$, and $\lambda_\sigma$ is a primitive $n$-th root of unity we say that $\sigma$ is \emph{purely nonsymplectic of order} $n$. 

The map $\Aut(X)\to \CC^*$ defined by $\sigma\mapsto \lambda_\sigma$ defines a character of $\Aut(X)$. For projective K3 surfaces, the image of this map is cyclic and we get the following short exact sequence (see e.g. \cite[Chapter 15]{huybrechtsK3}): 
\begin{equation}\label{eq:automorphism-ses-1}
    1\to \Aut(X)_{sym}\to \Aut(X) \to \mu_M\to 1 
\end{equation}
where $\Aut(X)_{sym}$ denotes the subgroup of symplectic automorphisms. 

Given any subgroup $G\subseteq \Aut(X)$, we get can restrict \eqref{eq:automorphism-ses-1} to the short exact sequence
\begin{equation}\label{eq:automorphism-short-exact-sequence}
    1\to G_{sym}\to G \to \mu_m\to 1 
\end{equation}
where $G_{sym}$ acts symplectically, and $\mu_m$ denotes the cyclic group of order $m$.

\subsection{K3 surfaces in weighted projective space}\label{ss:K3construction}

In this article we focus on K3 surfaces constructed from hypersurfaces in certain weighted projective spaces. In this section, we describe the polynomials and symmetry groups needed to construct such a K3 surface, and then we describe how to construct a K3 surface from these data. 

\subsubsection{Defining polynomials and symmetry groups}
We require a polynomial $W$ and a group of symmetries $G$ satisfying certain conditions, which we now describe. For more details on the construction and conditions, see e.g. \cite{ABS}, \cite{BCP}, \cite{ChR2}, \cite{CLPS}, \cite{CP}, \cite{krawitz}.


\begin{definition}
    \mbox{} 
    \begin{enumerate}
        \item A polynomial $W\in \mathbb{C}[x_1,\dots,x_n]$ is \emph{quasihomogeneous of degree $d$} if there exist positive integers $w_1, \cdots, w_n$ and $d$ such that
        \[
            W(c^{w_1}x_1, \cdots , c^{w_n}x_n) = c^d W(x_1, \cdots , x_n)
        \]
        for any $c\in\CC^*$. The tuple $(w_1, \cdots, w_n; d)$ is called the \emph{weight system} for $W$. We require that $\gcd(w_1,\dots,w_n)=1$.
        
        \item A quasihomogeneous polynomial is \emph{nondegenerate} if it has an isolated singularity at the origin and the weights are uniquely determined. A nondegenerate quasihomogeneous polynomial is called \emph{invertible} if it has the same number of monomials as variables. \label{quasihom-poly}
    \end{enumerate}
\end{definition}

It is a well-known fact (see e.g. \cite{KrSk}) that a polynomial is invertible if and only if it can be written as a Thom--Sebastiani sum of the following atomic types. 
\begin{itemize}
    \item Fermat: $x^a$ $\label{Fermat-general}$
    \item Chain: $x_1^{a_1}x_2+x_2^{a_2}x_3+\dots + x_n^{a_n}$  $\label{Chain-general}$
    \item Loop: $x_1^{a_1}x_2 + x_2^{a_2}x_3 + \cdots + x_n^{a_n}x_1$ $\label{Loops}$
\end{itemize}
We require $n>1$ for a chain or loop, and $a_i\geq 2$ for all exponents.  

Up to a change of variables, the exponents on each monomial determine an invertible polynomial. It is often useful to collect the exponents into a matrix $A_W$ which we call the \emph{exponent matrix}. The convention is that rows of $A_W$ represent monomials, and columns represent the variables. For example, the matrix
\[ 
\begin{pmatrix}
    a & 1 & 0\\
    0 & b & 1\\
    1 & 0 & c
\end{pmatrix}
\]
represents the (loop) polynomial $x^ay+y^bz+z^cx$.

The second ingredient is a group of symmetries $G$. Let $W$ be any nondegenerate quasihomogeneous polynomial. We define the maximal symmetry group
\[
 \Gmax_W = \{g \in \GL_n(\mathbb{C}) ~|~ W 
     (g\vec{x})=W(\vec{x})\text{ and } g_{ij}=0 \text{ unless }w_i=w_j\}.
\]
We will drop the subscript, when $W$ is clear from the context. There is an important subgroup of $\Gmax$, namely the group of diagonal symmetries, defined by
\[
    \Gdiag_W:=\{(g_1,\dots, g_n)\in (\CC^*)^n \mid W(g_1 x_1,\ldots,g_n x_n)=W(x_1,\ldots, x_n)\}. 
\]
Again, if $W$ is clear from the context, we will drop the subscript. One can view these diagonal symmetries as diagonal matrices, giving us an inclusion $\Gdiag \subseteq \Gmax$.

When working with diagonal symmetries, the following (additive) notation is convenient. For an invertible polynomial $W$, if $(g_1,g_2,\dots,g_n)\in \Gdiag$, then $g_j$ is a root of unity. So we can write $g_j=e^{2\pi i \Theta_j}$ for $\Theta_j\in \QQ/\ZZ$. We can think of $\Gdiag$ as a subgroup of $(\QQ/\ZZ)^n$ via the correspondence 
\[
(g_1,\dots, g_n)\leftrightarrow (\Theta_1,\dots,\Theta_n).
\]
Furthermore, the group $\Gdiag$ is generated (additively) by the columns of $A_W^{-1}$ (see e.g. \cite{krawitz}). Indeed, $|\Gdiag|=|\det A_W|$. 

There is a special element of $\Gdiag$ called the \emph{exponential grading operator}, which we denote (additively) by 
\[
j_W=(\tfrac{w_1}{d}, \dots, \tfrac{w_n}{d}).
\]
We denote the cyclic group of order $d$ generated by $j_W$ by $J_W := \langle j_W \rangle$. We note that $j_W$ is equal to the sum of the columns of $A_W^{-1}$. The condition $g_{ij} = 0$ unless $w_i = w_j$ in the definition of $\Gmax$ ensures $j_W$ is in the center of $\Gmax$ (see e.g. \cite{mukai}).

One more important subgroup of $\Gmax_W$ is $\SL_W:=\Gdiag\cap \SL(n,\CC)$, which we will encounter again later.

\subsubsection{Constructing a K3 surface $X_{W,G}$}

Reid (unpublished) and Yonemura in \cite{yonemura} classified 95 weight systems for which a nondegenerate polynomial that is quasihomogeneous with respect to one of these weight systems yields a K3 surface. 
We now describe how to construct a K3 surface given a pair $(W,G)$ as above, where $W\in \CC[x,y,z,w]$ is an invertible polynomial that is quasihomogeneous with respect to one of the 95 weight systems and $G$ is a group of diagonal symmetries satisfying $J_W \subseteq G\subseteq \SL_W$. 

Let $V(W)$ denote the vanishing locus of $W$. For a general quasihomogeneous polynomial, $V(W)$ will be a singular surface. However, since $W$ is nondegenerate, the singularities are all inherited from the ambient space. We can define the K3 surface $X_W=\widetilde{V(W)}$ as the minimal resolution of $V(W)$, but we also want to incorporate $G$ into our construction. 

To this end, notice the group $\overline{G}=G/J_W$ acts naturally on $V(W)$, extending to an action on $X_W$. Since $G\subseteq \SL_W$, the group $\overline{G}$ acts symplectically on $X_W$, preserving the volume form (see \cite[Proposition~2.3]{ABS}).
Again the minimal resolution
\[
X_{W,G}:=\widetilde{X_W/\overline{G}}
\]
is a K3 surface. Such K3 surfaces are the objects of our study.

\subsection{Lattice Theory Background}\label{ss:latticebackground}

Before we can introduce mirror symmetry we provide some background with integral lattices. This discussion follows \cite[Section 1.2]{CP}. For more details, the interested reader should consult Nikulin's groundbreaking paper \cite{nikulin}.

We begin with some preliminary definitions. An \emph{integral lattice}, or lattice for short, is a free abelian group $L$ of finite rank together with a non-degenerate, symmetric (integral) bilinear form $B:L\times L\rightarrow \mathbb{Z}$. A lattice is \emph{even} if the associated bilinear form satisfies $B(x,x)\in 2\ZZ$ for all $x\in L$. The \emph{signature} $(t_+,t_-)$ of a lattice, is the signature of $B$.
We say $L\subseteq L'$ is a \emph{primitive sublattice} of $L'$ if $L'/L$ is free. On the other hand, we say $L'\supseteq L$ is an \emph{overlattice} of $L$ (of finite index) if $L'/L$ is a finite abelian group. 
Given an even lattice, we can embed $L\hookrightarrow L^\vee:=\Hom(L,\ZZ)$ via the bilinear form $B$. We define the \emph{discriminant group} $A_L:=L^\vee/L$, which is a finite group of order $|\det B|$. An even lattice is \emph{unimodular} if $A_L$ is trivial. The \emph{length} $l(A_L)$ is the minimum number of generators of $A_L$. 

The bilinear form $B$ induces a quadratic form $Q:L\to 2\ZZ$ by $Q(x) = B(x, x)$. This quadratic form descends to the \emph{discriminant quadratic form}
\[
q_L:A_L\to \QQ/2\ZZ 
\]
giving us a finite quadratic form.

The set of finite quadratic forms has the structure of a semi-group via the operation of direct sum. The semi-group is generated by the following forms, as proven in \cite{nikulin} (see also \cite{belcastro} for a nice treatment):
\begin{enumerate}
    \item When $p\neq 2$ is prime, with $k\in\ZZ^+$, and $\epsilon\in\{\pm 1\}$, we define the quadratic form 
    \[
    \omega^\epsilon_{p,k}:\ZZ/p^k\ZZ\rightarrow \QQ/2\ZZ
    \]
    by $\omega^\epsilon_{p,k}(1)=ap^{-k}$ where $a$ is the smallest even integer with quadratic residue $\epsilon$ modulo $p$.
    \item When $p=2$, $k\in\ZZ^+$, and $\epsilon\in\{\pm 1,\pm 5\},$ we define 
    \[
    \omega^\epsilon_{2,k}:\ZZ/2^k\ZZ\rightarrow \QQ/2\ZZ
    \]
    by $\omega^\epsilon_{2,k}(1)=\epsilon\cdot 2^{-k}$.
    \item For $k\in\ZZ^+$, we define $u_k$ and $v_k$ on $\ZZ/2^k\ZZ\times \ZZ/2^k\ZZ$ by the matrices
    \begin{align*}
        u_k=\begin{pmatrix}
            0&2^{-k}\\2^{-k}&0
        \end{pmatrix},\quad v_k=\begin{pmatrix}
            2^{-k+1}&2^{-k}\\2^{-k}&2^{-k+1}
        \end{pmatrix}.
    \end{align*}
\end{enumerate}
Every finite quadratic form may be written (not uniquely) as the direct sum of the generators $\omega^\epsilon_{p, k}, u_k, v_k$.

Given a prime $p$, we say a lattice $L$ is \emph{$p$-elementary} if the discriminant group satisfies $A_L\cong (\ZZ/p\ZZ)^a$ for some $a>0$. In this case, $a$ is the length of $L$.

We now describe several important lattices. 
First, the lattice $U$ is the unimodular, indefinite lattice of rank 2, whose symmetric form can be represented by the matrix 
\[
\begin{pmatrix}
    0 & 1\\ 1 & 0
\end{pmatrix}.
\]
The lattices $A_n$ for $n\geq 2$, $D_n$ for $n\geq 4$, and $E_6, E_7, E_8$ are negative definite lattices, whose symmetric form can be represented by the corresponding Dynkin diagram. In other words, the corresponding lattice has a basis $e_1,e_2, \dots, e_n$ corresponding to the vertices of the diagram, with $B(e_i,e_i)=-2$, and $B(e_i,e_j)=1$ if the $i$-th and $j$-th vertices share an edge. $B(e_i,e_j)$ is zero otherwise. 

For example, $E_8$ has Dynkin diagram 
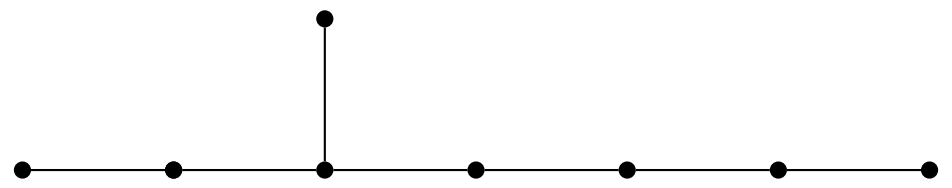
\begin{center}
\begin{tikzpicture}[
    every node/.style={
        circle,
        draw,
        fill=black,
        inner sep=1pt,
    }
]

\foreach \x/\name in {
    0/a,
    1/b,
    2/c,
    3/d,
    4/e,
    5/f,
    6/g,
} {
    \node (\name) at (\x,0) {};
}

\node (h) at (2,1) {};

\draw (a)--(b)--(c)--(d)--(e)--(f)--(g);
\draw (c)--(h);

\end{tikzpicture}
    
\end{center}
so the corresponding lattice is negative definite with rank 8 and symmetric bilinear form given by the diagram, as described above. 

There are three important lattices related to any K3 surface that will be relevant in this article. First, the cohomology group $H^2(X,\ZZ)$ is an even unimodular lattice of signature $(3,19)$. As such $H^2(X,\ZZ)\cong U^3\oplus (E_8)^2$. This is sometimes referred to as the \emph{K3 lattice}. 

If $X$ is an algebraic K3 surface, then the \emph{Picard lattice} $\Pic(X)$ is an even lattice of signature $(1,t)$ for some nonnegative integer $t$. Indeed $\Pic(X)$ is a primitive sublattice of $H^2(X,\ZZ)$. 

Furthermore, if $X$ has a nonsymplectic automorphism of order $n$, then the invariant lattice $S(\sigma):=\{x\in H^2(X,\ZZ) \mid \sigma^*x=x\}$ is a primitive sublattice of $\Pic(X)$. As we have seen, the invariant lattice seems to be the key to understanding LPK3 mirror symmetry. We will discuss this more in the next section. 

Again, let $p$ be a prime, and suppose $\sigma$ is a purely nonsymplectic automorphism of $X$ of order $p$. It turns out that the invariant lattice $S(\sigma)$ is $p$-elementary. In particular, if $\sigma$ is a nonsymplectic involution, then $S(\sigma)$ is 2-elementary (see e.g. \cite{MO} and \cite{other_primes}). There are several nonsymplectic involutions, which will play an important role in our main result.

There is a classification of 2-elementary lattices due to Nikulin and there are similar results for any prime $p$ (see e.g. \cite{other_primes}).
\begin{theorem}[cf. Nikulin \cite{nikulin2}]\label{thm:2lattice}
    An even, indefinite, 2-elementary lattice is determined by the rank, length, and a third invariant $\delta\in \{0,1\}$. 
\end{theorem}

\begin{figure}
    \centering
\begin{tikzpicture}[
    x=0.48cm,
    y=0.48cm,
    every node/.style={font=\normalsize},
  ]

  \draw[->,line width=0.35pt] (0,0) -- (22.25,0);
  \draw[->,line width=0.35pt] (0,0) -- (0,12.05);

  \foreach \r in {1,...,20} {
    \draw[line width=0.3pt] (\r,-0.10) -- (\r,0.10);
    \node[below=3pt] at (\r,0) {$\r$};
  }
  \foreach \a in {1,...,11} {
    \draw[line width=0.3pt] (-0.10,\a) -- (0.10,\a);
    \node[left=4pt] at (0,\a) {$\a$};
  }
  \node[below left=3pt] at (0,0) {$0$};
  \node[above left=2pt] at (0,12.05) {$a$};
  \node[below right=2pt] at (22.25,0) {$r$};

  \draw[->,line width=0.65pt]
    (11,11) -- (0.05,0.05);
  \draw[->,line width=0.65pt]
    (11,11) -- (21.55,0.45);

  \node at (11,11.58) {$0$};
  \foreach \n/\r in {1/10,2/9,3/8,4/7,5/6,6/5,7/4,8/3,9/2,10/1} {
    \node[anchor=east] at (\r,\r+0.42) {$\n$};
  }
  \foreach \n/\r in {1/12,2/13,3/14,4/15,5/16,6/17,7/18,8/19,9/20} {
    \node[anchor=west] at (\r+0.23,22-\r+0.42) {$\n$};
  }

  \node[anchor=west] at (-0.1,0.62) {$g$};
  \node[anchor=west] at (21.43,0.67) {$k$};

  \foreach \r/\a in {
      1/1,3/3,4/4,5/5,6/6,7/7,8/8,9/9,10/10, 11/11,
      12/10,13/9,14/8,15/7,16/6,17/5,19/3,20/2
    } {
    \fill (\r,\a) circle (1.65pt);
  }

  \foreach \r/\a in {
      3/1,9/1,11/1,17/1,19/1,
      4/2,8/2,12/2,16/2,
      5/3,7/3,9/3,11/3,13/3,15/3,17/3,
      8/4,12/4,16/4,
      7/5,9/5,11/5,13/5, 15/5,
      8/6,10/6, 12/6, 14/6, 
      9/7,11/7,13/7, 10/8, 12/8,
      11/9, 2/2,10/2,18/2,6/4,10/4,14/4,18/4
    } {
    \fill (\r,\a) circle (1.65pt);
  }

  \foreach \r/\a in {
      2/1.8,6/2,10/1.8,14/2,18/1.8,
      6/3.8,10/3.8,14/3.8,18/3.8,
      10/5.8,14/5.8,
      10/7.8,10/9.8, 10/0, 2/0, 18/0
    } {
    \node[inner sep=0pt] at (\r,\a) {$\ast$};
  }

  \fill (18.85,12.02) circle (1.65pt);
  \node[anchor=west] at (19.55,12.02) {$\delta=1$};
  \node[inner sep=0pt] at (18.85,11.05) {$\ast$};
  \node[anchor=west] at (19.55,11.05) {$\delta=0$};

\end{tikzpicture}
\caption{Nikulin's classification of invariant lattices for nonsymplectic involutions.}
\label{fig:2elementary}
\end{figure}
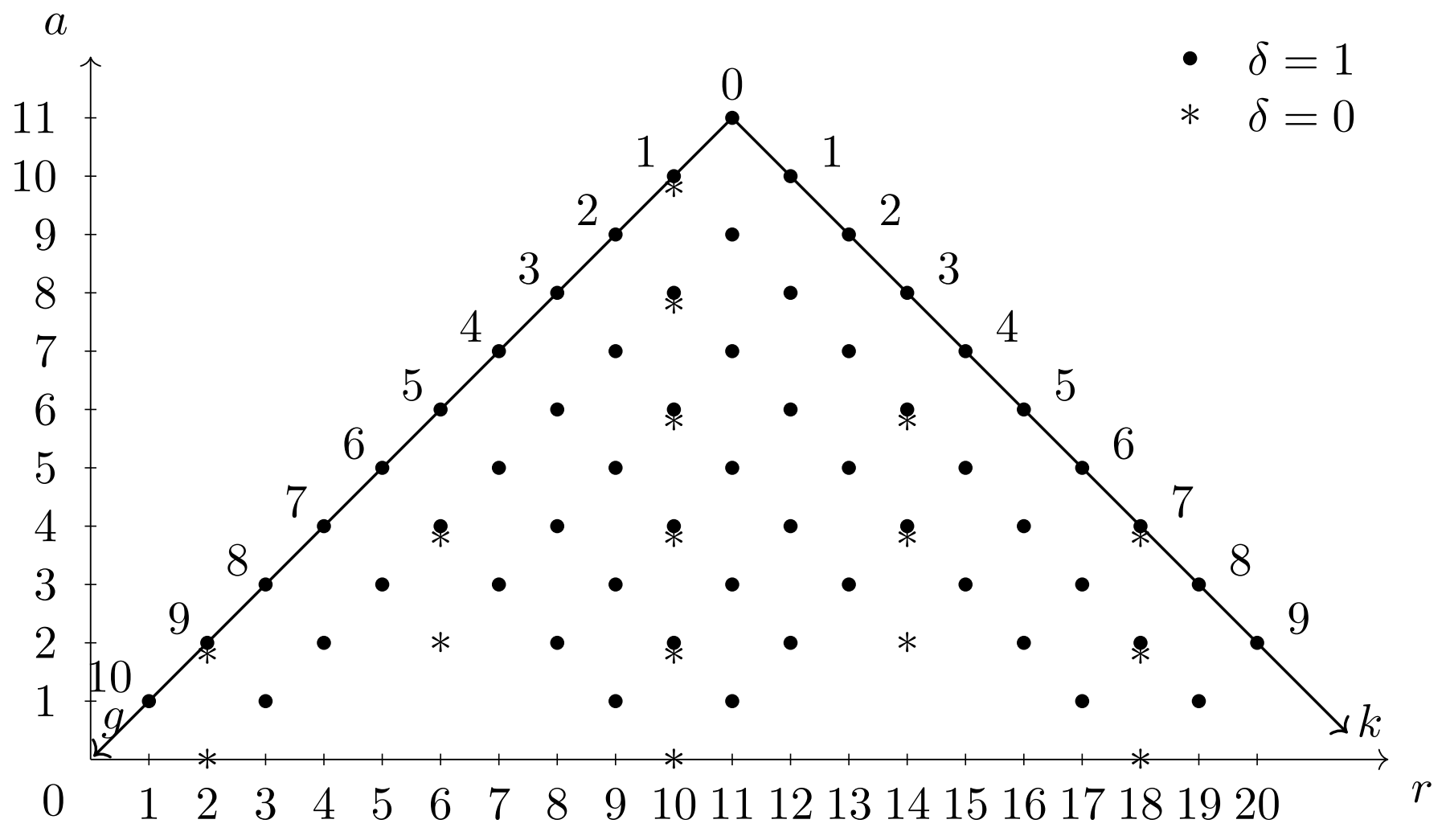

In defining $\delta$, we say $\delta=0$ if $q_L(x)\in \ZZ$ for every element $x\in A_L$, otherwise $\delta=1$. There is a classification of even, indefinite, 2-elementary lattices given in \cite{nikulin4} and \cite{nikulin2}, which we have reproduced in Figure~\ref{fig:2elementary} (see also \cite{other_primes}). Here $r$ is the rank, and $a$ is the length, i.e. $A_L\cong (\ZZ/2\ZZ)^a$. 

Furthermore, the invariant lattice of a nonsymplectic involution can be completely determined by the topology of its fixed locus, as described in the following theorem of Nikulin: 
\begin{theorem}[cf. Nikulin \cite{nikulin2}]\label{thm:nikulingark}
    The fixed locus of a nonsymplectic involution on a K3 surface is one of the following:
    \begin{itemize}
        \item empty if $r=10$, $a=10$ and $\delta=0$,
        \item the disjoint union of two elliptic curves if $r=10$, $a=8$, and $\delta=0$,
        \item the disjoint union of a curve of genus $g$ and $k$ rational curves otherwise, where $g=(22-r-a)/2$, $k=(r-a)/2$. 
    \end{itemize}
\end{theorem}
Thus, one way to completely determine the invariant lattice of a nonsymplectic involution is by finding the topology of the fixed locus, and this is also depicted in Figure~\ref{fig:2elementary}. We will leverage this fact later. 

\subsection{Mirror Symmetry} \label{sec:mirror-symmetry}

Several constructions of mirror symmetry have been developed in the literature, each arising from a different mathematical context. The K3 surfaces in this article naturally fit two such constructions for which notions of mirror symmetry have been defined, namely BHK mirror symmetry and LPK3 mirror symmetry.
As mentioned previously, a natural question one might pose is whether the two versions of mirror symmetry produce the same mirror. Before we address this question, we describe each version of mirror symmetry. 

\subsubsection{BHK mirror symmetry}\label{ss:bhk}
Although we have described mirror symmetry in geometric terms, BHK mirror symmetry was first described in the setting of what are known as Landau--Ginzburg (LG) models. Much has been written on the topic (see e.g. \cite{ChK}, \cite{CJMPW}, \cite{FJJS}, \cite{FPS}, \cite{krawitz}, \cite{KP}), so we will be brief. 

In this setting, the input data for an LG model is a pair $(W,G)$, where $W$ is a nondegenerate quasihomogeneous polynomial, and $G\subseteq G^{diag}$ is a group of symmetries of $W$. Given such a pair, if $W$ is invertible, the BHK mirror symmetry construction returns another pair $(W^T,G^T)$, defined as follows. 

The polynomial $W^T$ is obtained by exchanging the roles of monomials and variables. More specifically, given an invertible polynomial
\[
W(x_0, \ldots, x_n) = \sum_{i=0}^n \prod_{j=0}^n x_j^{a_{ij}},
\]
we can encode the exponents in the exponent matrix $A_W$ as described in Section~\ref{ss:K3construction}. We define the dual polynomial $W^T$ to be the invertible polynomial whose exponent matrix is obtained from the transpose $(A_W)^T$. Recall every invertible polynomial is a sum of atomic types (see Section~\ref{ss:K3construction}). One might observe that this construction preserves each atomic type, but reverses the direction of loops and chains. 

For example, we have seen that the loop polynomial $x^ay+y^bz+z^cx$ has matrix 
\[
A_W=\begin{pmatrix}
    a & 1 & 0\\
    0 & b & 1\\
    1 & 0 & c
\end{pmatrix}
\]
Then we have 
\[
A_W^T=\begin{pmatrix}
    a & 0 & 1\\
    1 & b & 0\\
    0 & 1 & c
\end{pmatrix}
\]
which represents the polynomial $x^az+xy^b+yz^c$.

The dual group $G$ is defined (additively) as
\[
G^T = \{g \in \Gdiag_{W^T}: gA_Wh^T \in \ZZ \quad \text{for all} \quad h \in G\}.
\]
In practice, this definition is difficult to work with. However, we can often determine the dual group from the following properties, which are proved in \cite{ABS}:
\begin{enumerate}
    \item $J_W^T=\SL_{W^T}$ \label{prop:jdual}
    \item $(\Gdiag_W)^T$ is the trivial group \label{prop:dualdual}
    \item $(G^T)^T=G$ \label{prop:maxdual}
    \item If $G_1\subseteq G_2$, then $G_2^T\subseteq G_1^T$ and $G_1^T/G_2^T\cong G_2/G_1$. \label{prop:quotientdual}
\end{enumerate}
 
In particular, if we begin with a group $G$ satisfying
\[
J_W \leq G \leq SL_W
\]
then these properties also give us 
\[
J_{W^T} \leq G^T \leq SL_{W^T}.
\]

If we start with an invertible polynomial $W$, quasihomogeneous with respect to one of the 95 weight systems, and a group $G$ satisfying $J_W\subseteq G\subseteq \SL_W$, then remarkably the dual pair also satisfies the same relationships, namely $W^T$ is invertible with respect to one of the 95 weight systems, and $J_{W^T} \leq G^T \leq SL_{W^T}$. Thus we obtain another K3 surface $X_{W^T,G^T}$. 

If $(W,G)$ satisfy certain conditions, then $X_{W,G}$ is a Calabi--Yau (CY) variety. 
There is another deep relationship between the LG model $(W,G)$, and the CY manifold $X_{W,G}$, known as the LG/CY correspondence, which relates the LG model for a pair $(W,G)$ with the enumerative geometry of the CY manifold $X_{W,G}$ (see e.g. \cite{ChIR}, \cite{ChR1}, \cite{ChR3}, \cite{ChR2}, \cite{LPSh}, \cite{PSh}).
Although BHK originated as an explicit algebraic duality on the LG models, when we combine it with  the LG/CY correspondence, we obtain a mirror construction for Calabi--Yau manifolds of the form $X_{W,G}$, namely BHK mirror symmetry says that the K3 surface $X_{W,G}$ is mirror dual to the K3 surface $X_{W^T,G^T}$.

\subsubsection{LPK3}
The second form of mirror symmetry we consider is called Lattice-polarized K3 (LPK3) mirror symmetry. This was originally described by Dolgachev in \cite{dolgachev}. This version of mirror symmetry formulates mirror symmetry intrinsically in terms of the K3 lattice
\begin{align*}
    L_{K3} &= U^3 \oplus (E_8)^2.
\end{align*}
The idea of this version of mirror symmetry came from the classical rotation of the Hodge diamond that is well-known for Calabi--Yau threefolds. However, rotating the Hodge diamond \eqref{eq:K3HodgeDiamond} for a K3 surface always leaves the Hodge diamond invariant. So the idea is to break up $H^{1,1}(X)$ into two parts that would be exchanged under mirror symmetry, much like how $H^{2,1}(X)$ and $H^{1,1}(X)$ are for a CY three-fold. 

In order to describe this version of mirror symmetry, we say a K3 surface $X$ is $M$-polarized if there exists a primitive embedding $j: M \to \text{Pic}(X)$ of a lattice $M$ of signature $(1,t)$. Its mirror lattice $M^V$ is defined by the orthogonal decomposition of the orthogonal complement
\[
M^{\perp}_{L_{K3}} \cong U \oplus M^V.
\]
if such a decomposition exists. 
Thus, we can define two K3 surfaces $X$ and $X^\vee$ to be LPK3 mirrors if $X$ and $X^\vee$ are polarized respectively by $M$ and $M^V$. In this way, we have ``split'' $H^{1,1}(X)$ into two parts, and mirror symmetry can be viewed as some sort of exchange between $M$ and $M^\vee$. 

By considering the rank of the lattices involved in the previous discussion, we see that a necessary condition for $M$ and $M^\vee$ to be mirror lattices is 
\begin{equation}\label{eq:mirrorRank}
\rk M^\vee =20-\rk M.    
\end{equation}

In practice, the polarization $M$ is often taken as the invariant lattice 
\[
S(\sigma) = \{x \in H^2(X,\ZZ):\sigma^* x = x\}
\]
of some purely nonsymplectic automorphism $\sigma$ of $X$. The mirror lattice then arises as the invariant lattice of a particular purely nonsymplectic automorphism $\sigma^T$ on the BHK mirror surface. 

To give sufficient conditions for two lattices to be mirror lattices is much more difficult, and a full discussion is not necessary for the current article. However, when $\sigma$ is a nonsymplectic involution, the situation becomes much simpler. In that case, recall that $S(\sigma)$ is a 2-elementary lattice, and is therefore determined (up to isomorphism) by the invariants $r, a, \delta$, where $r$ is the rank, $a$ is the length, and $\delta$ is as in Theorem \ref{thm:2lattice}. In that case, the mirror lattice is the (unique) 2-elementary lattice with invarants $20-r$, $a$, and $\delta$, if such a lattice exists (see Figure~\ref{fig:2elementary}, also \cite{nikulin2} and \cite{ABS}). We note that if we draw a line at $r=10$ in Figure~\ref{fig:2elementary},  then the mirror lattice is the lattice we obtain by reflecting across this line, with the same $\delta$ invariant.

\section{An Example} \label{sec:example} 
\subsection{Overview} 

We now turn our attention to an example of a sum of two chain potentials to show that Conjecture~\ref{conj:main} is false for diagonal symmetries. We will see that no diagonal symmetry will provide a suitable polarizing lattice. However, we will show that a nondiagonal nonsymplectic automorphism does.  

\subsubsection{Polynomial and Exponent Matrix}

Consider the polynomial
\[
    W=x^3y+xy^3+z^3w+zw^3~. \label{eq:definition-W}
\]
This polynomial is invertible. In fact, it is smooth and homogeneous. The exponent matrix is
\[
    A_W =
    \begin{pmatrix}
        3 & 1 & 0 & 0 \\
        1 & 3 & 0 & 0 \\
        0 & 0 & 3 & 1 \\
        0 & 0 & 1 & 3
    \end{pmatrix}
    ~.\label{eq:example-exponential-matrix}
\]
Since $A_W^T = A_W$, the polynomial $W$ is self-dual. Of course, for BHK mirror symmetry, we also need to consider a group of symmetries. As described in Section~\ref{ss:K3construction}, this means we need to understand the groups $G$ that satisfy $J_W\subseteq G\subseteq \SL_W$.

\subsubsection{Diagonal Symmetry Group}

With $W$ as in \eqref{eq:definition-W}, we have 
\[A_W^{-1} = 
    \begin{pmatrix}
        \frac{3}{8} & -\frac{1}{8} & 0 & 0 \\
        -\frac{1}{8} & \frac{3}{8} & 0 & 0 \\
        0 & 0 & \frac{3}{8} & -\frac{1}{8} \\
        0 & 0 & -\frac{1}{8} & \frac{3}{8}
    \end{pmatrix}\]
Recall $\Gdiag$ is generated by columns of $A_{W}^{-1}$, so
\[\Gdiag = \left\langle \left(\frac{3}{8}, -\frac{1}{8}, 0, 0\right), \left(-\frac{1}{8}, \frac{3}{8}, 0, 0\right), \left(0, 0, \frac{3}{8}, -\frac{1}{8}\right), \left(0, 0, -\frac{1}{8}, \frac{3}{8}\right)\right\rangle~.\]
We also have $j_W= \left(\frac{1}{4}, \frac{1}{4}, \frac{1}{4}, \frac{1}{4}\right)$.

Notice that
\begin{align*}
    5\left(\frac{3}{8}, -\frac{1}{8}, 0, 0\right) = \left(-\frac{1}{8}, \frac{3}{8}, 0, 0\right) \\
    5\left(0, 0, \frac{3}{8}, -\frac{1}{8}\right) = \left(0, 0, -\frac{1}{8}, \frac{3}{8}\right).
\end{align*}
So, by defining
\begin{align}
    \sigma_L &:= \left(\frac{3}{8}, -\frac{1}{8}, 0, 0\right) \label{eq:definition-sigmaL} \\
    \sigma_R &:= \left(0, 0, \frac{3}{8}, -\frac{1}{8}\right), \label{eq:definition-sigmaR}
\end{align}
we find
\[
\Gdiag = \left\langle \sigma_L, \sigma_R\right\rangle.
\]
Since each generator is order 8 and they are linearly independent, we conclude
\[
\Gdiag \cong \ZZ_8 \oplus \ZZ_8,
\]
and $\{\sigma_L, \sigma_R\}$ is a minimal set of generators. It follows $|\Gdiag| = 64$. 

\subsubsection{Some nondiagonal symmetries}\label{ss:nondiagonal}

Although Conjecture~\ref{conj:main} was originally formulated for diagonal symmetries, there are other natural nondiagonal symmetries that we can also consider which come from permuting the variables of $W$. These will also lead to automorphisms of our K3 surfaces. 

To describe the subgroup of $S_4$ that arises from permuting the variables while also preserving $W$, we will use cycle notation to denote the symmetries, e.g. $(xy)$ will denote the symmetry that exchanges the variable $x$ for the variable $y$. Geometrically, this exchanges the $x$ coordinate and $y$ coordinate of any point of $V(W)$. 

There are two kinds of generators for the group in question---those permuting the variables within a loop---in our example this would be $(xy)$ and $(zw)$---, and the one that exchanges the two loops---this would be $(xz)(yw)$. As shown in \cite{CJMPW}, the largest group of permutations for $W$ is the wreath product $(\ZZ/2\ZZ) \wr (\ZZ/2\ZZ)$. Explicitly, there are eight permutations in this group 
\[
\{(1), (xy), (zw), (xy)(zw), (xz)(yw), (xzyw), (xwyz), (xw)(yz)\}.
\]

Each of these automorphisms will induce an automorphism on the K3 surfaces we consider. One can check that such an automorphism is symplectic if and only if the corresponding permutation is even. 

\subsection{Three K3 surfaces}

Let us next examine the group $\SL_W$. Recall that each K3 surface we consider is constructed taking the quotient of $V(W)$ by a group $\overline{G}$, with $G \subseteq \SL_W$. Naturally, one may ask which groups $\overline{G}$ are possible. We can use properties introduced earlier to answer this question.

From Properties~\eqref{prop:jdual}--\eqref{prop:quotientdual} in Section~\ref{ss:bhk}, we have, for example, $\Gdiag/\SL_W\cong J_{W^T}$. Since $W^T=W$, we also have $J_{W^T}=J_W$. Hence we get the following lattice of subgroups: 
\[
\begin{array}{cccc}
     &\Gdiag & 1 & \\
      & | & | & 4\\
     & \SL_W & J_W& \\
      & | & | & \\
     & J_W & \SL_W & \\
     4 & | & | & \\
      & 1 & \Gdiag & 
\end{array}
\]
The lattice on the right is inverted to indicate the index of subgroups going up is isomorphic to the index going down on the left. We can fill in the remaining group indices quite simply. First we have $|\Gdiag/\SL_W|=|J_{W^T}|=4$, as previously mentioned. Then we can also fill in $|\SL_W/J_W|=4$, since $|\Gdiag|=4^3$.

Notice $\sigma_L-\sigma_R=\left(\frac{3}{8},-\frac{1}{8},-\frac{3}{8},\frac{1}{8}\right)$ has order 4 in $\SL_W/J_W$, so $\SL_W/J_W\cong \ZZ/4\ZZ$. Hence there are exactly three subgroups $G$ satisfying $J_W\subseteq G\subseteq\SL_W$. These are 
\begin{align}
    G_1&=J_W =\left\langle\left(\frac{1}{4},\frac{1}{4},\frac{1}{4},\frac{1}{4}\right)\right\rangle~, \nonumber\\
    G_2&=\left\langle\left(\frac{1}{4},\frac{1}{4},\frac{1}{4},\frac{1}{4}\right), \left(\frac{1}{2},\frac{1}{2},0, 0\right)\right\rangle~, \label{eq:definition-of-G2} \\
    G_4&=\SL_W =\left\langle\left(\frac{1}{4},\frac{1}{4},\frac{1}{4},\frac{1}{4}\right), \left(\frac{3}{8},-\frac{1}{8},-\frac{3}{8},\frac{1}{8}\right)\right\rangle~. \nonumber 
\end{align}
The subscript on the group indicates the order of $\overline{G}_i$, e.g. $G_2/J_W$ has order 2. 
Given $W$ introduced above, for each $G_i$ we can now construct the associated K3 surface, which we will denote
\[
X_i := X_{W, G_i}~.
\]

We have seen that $W^T=W$, and again from Properties~\eqref{prop:jdual}--\eqref{prop:quotientdual} in Section~\ref{ss:bhk}, one can easily work out that $G_1^T = G_4$ and $G_2^T=G_2$. So according to BHK mirror symmetry, we have two pairs of mirror surfaces, namely $X_1$ and $X_4$ are BHK mirror duals, and $X_2$ is mirror dual to itself.

In the following sections, we will first describe $X_1$ and $X_4$, and the invariant lattice for any nonsymplectic diagonal symmetry, then we will do the same with $X_2$. 

\subsection{Surface $X_1$}

\subsubsection{Geometry}

As before, we let $V(\cdot)$ denote the vanishing set in $\CP^3$. For $W$ introduced in \eqref{eq:definition-W}, notice $V(W) \subseteq \CP^3$ has no singularities. Furthermore, $\overline{G}_1$ is trivial, so $X_1=V(W)$. 

We find a set of lines which form a basis for the Picard lattice. It turns out that $X_1$ is isomorphic to the Fermat quartic $V(x^4 + y^4 + z^4 + w^4)\subseteq\CC\PP^3$, so the Picard lattice is generated by 48 lines (see e.g. \cite{SSvL}). To find the lines on $X_1$, we take the image of the 48 known lines on the Fermat quartic under the isomorphism (see Appendix~\ref{appendix:A}). 

We introduce some notation to simplify this exposition. Let $\zeta$ denote the primitive eighth root of unity and define the following polynomials
\begin{align*}
    \aplus &:= x + y & \bplus &:= z + w \\
    \aminus &:= x - y & \bminus &:= z - w~.
\end{align*}
There are 3 families of 16 lines each on $X_1$. For $m, n \in \ZZ/4\ZZ$ we denote these as
\begin{align*}
    C_1(m, n) &:= V(\aplus + i^m\aminus,~\bplus + i^n\bminus) \\
    C_2(m, n) &:= V(\aplus + i^m\zeta \bplus,~\aminus + i^n\zeta \bminus) \\
    C_3(m, n) &:= V(\aplus + i^m\bminus,~\aminus + i^n\bplus)~.
\end{align*}
It turns out the first family $C_1$ is particularly important for this example. So we will give these lines a more intuitive description: 
By defining the polynomials
\begin{align*}
    \alpha &= x+iy & \bar{\alpha} &= x - iy \\
    \beta &= z+iw & \bar{\beta} &= z - iw~,
\end{align*}
we can rewrite each line in family $C_1$ using subscripts to notate the defining linear forms. For example, we may write $C_{x\beta}:=V(x,\beta)$. Thus we have the following alternate description of family $C_1$:
\begin{align*}
    C_{xz} &= C_1(0, 0) & C_{yz} &= C_1(2, 0) & C_{\alpha z} &= C_1(3, 0) & C_{\alpha\beta} &= C_1(3, 3) \\
    C_{xw} &= C_1(0, 2) & C_{yw} &= C_1(2, 2) & C_{\bar{\alpha} z} &= C_1(1, 0) & C_{\alpha\bar{\beta}} &= C_1(3, 1) \\
    C_{x\beta} &= C_1(0, 3) & C_{y\beta} &= C_1(2, 3) & C_{\alpha w} &= C_1(3, 2) & C_{\bar{\alpha}\beta} &= C_1(1, 3)\\
    C_{x\bar{\beta}} &= C_1(0, 1) & C_{y\bar{\beta}} &= C_1(2, 1) & C_{\bar{\alpha} w} &= C_1(1, 2) & C_{\bar{\alpha}\bar{\beta}} &= C_1(1, 1)~.
\end{align*}

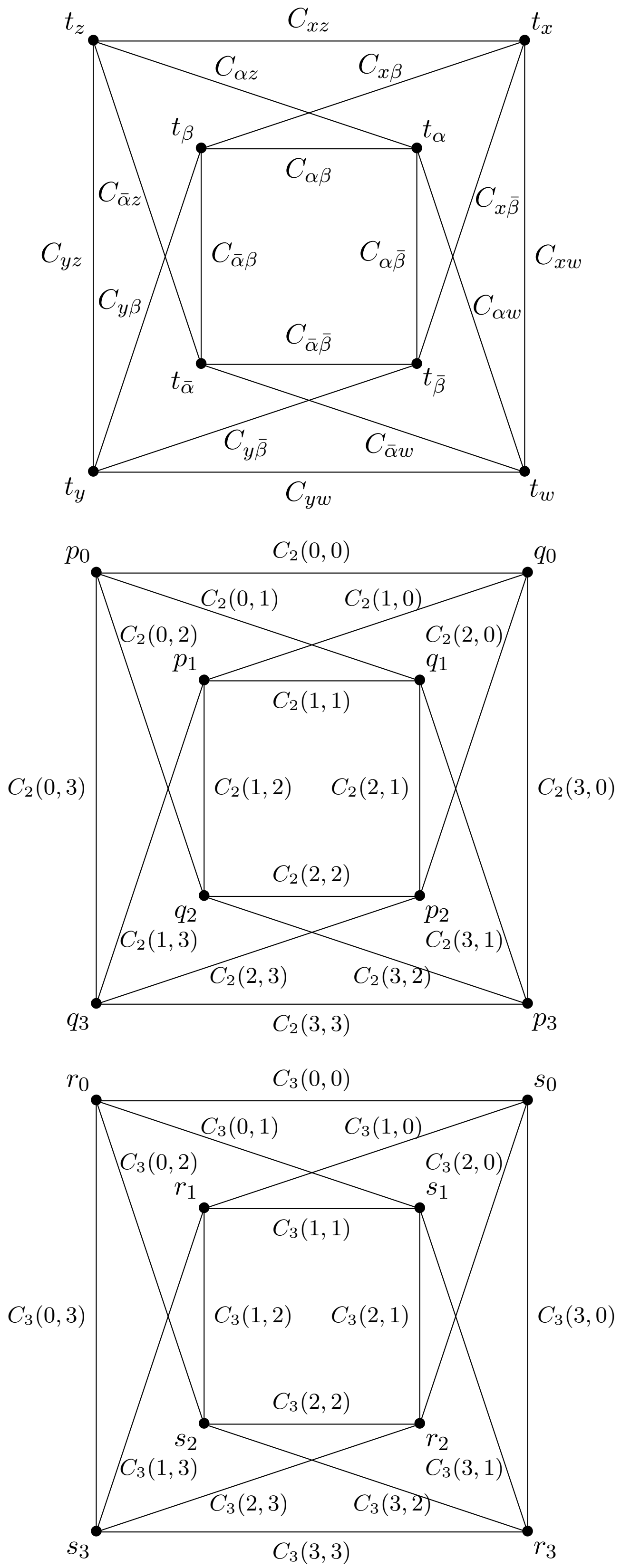
\begin{figure}
\centering\begin{tikzpicture}[xscale=.5,yscale=.5]
\draw (-6,6)--(6,6)--(6,-6)--(-6,-6)--(-6,6);
\foreach \Point in {(-6,6),(6,6),(6,-6),(-6,-6)}{
    \node at \Point {\textbullet};
}
\node [above] at (0,6){$C_{xz}$};
\node at (-6.5,6.5){$t_z$};
\node [right] at (6,0){$C_{xw}$};
\node at (6.5,6.5){$t_x$};
\node [below] at (0,-6){$C_{yw}$};
\node at (6.5,-6.5){$t_w$};
\node [left] at (-6,0){$C_{yz}$};
\node at (-6.5,-6.5){$t_y$};

\draw (-3,3)--(3,3)--(3,-3)--(-3,-3)--(-3,3);
\foreach \Point in {(-3,3),(3,3),(3,-3),(-3,-3)}{
    \node at \Point {\textbullet};
}
\node [below] at (0,3){$C_{\alpha\beta}$};
\node at (3.5, 3.5){$t_\alpha$};
\node [left] at (3,0){$C_{\alpha\bar{\beta}}$};
\node at (3.5, -3.5){$t_{\bar{\beta}}$};
\node [above] at (0,-3){$C_{\bar{\alpha}\bar{\beta}}$};
\node at (-3.5, -3.5){$t_{\bar{\alpha}}$};
\node [right] at (-3,0){$C_{\bar{\alpha}\beta}$};
\node at (-3.5, 3.5){$t_\beta$};

\draw (-6,6)--(3, 3);
\node at (-2, 5.25){$C_{\alpha z}$};
\draw (-6,6)--(-3, -3);
\node at (-5.25, 1.75){$C_{\bar{\alpha}z}$};
\draw (6,6)--(-3, 3);
\node at (2, 5.25){$C_{x\beta}$};
\draw (6,6)--(3, -3);
\node at (5.25, 1.5){$C_{x\bar{\beta}}$};
\draw (6,-6)--(3, 3);
\node at (5.25, -1.4){$C_{\alpha w}$};
\draw (6,-6)--(-3, -3);
\node at (2.25, -5.25){$C_{\bar{\alpha}w}$};
\draw (-6,-6)--(-3, 3);
\node at (-5.25, -1.3){$C_{y\beta}$};
\draw (-6,-6)--(3, -3);
\node at (-1.75, -5.25){$C_{y\bar{\beta}}$};
\end{tikzpicture}

\vspace{0.15cm}

\begin{tikzpicture}[xscale=.5,yscale=.5]
\draw (-6,6)--(6,6)--(6,-6)--(-6,-6)--(-6,6);
\foreach \Point in {(-6,6),(6,6),(6,-6),(-6,-6)}{
    \node at \Point {\textbullet};
}
\node [above] at (0,6){\footnotesize$C_2(0,0)$};
\node at (-6.5,6.5){$p_0$};
\node [right] at (6,0){\footnotesize$C_2(3,0)$};
\node at (6.5,6.5){$q_0$};
\node [below] at (0,-6){\footnotesize$C_2(3,3)$};
\node at (6.5,-6.5){$p_3$};
\node [left] at (-6,0){\footnotesize$C_2(0,3)$};
\node at (-6.5,-6.5){$q_3$};

\draw (-3,3)--(3,3)--(3,-3)--(-3,-3)--(-3,3);
\foreach \Point in {(-3,3),(3,3),(3,-3),(-3,-3)}{
    \node at \Point {\textbullet};
}
\node [below] at (0,3){\footnotesize$C_2(1,1)$};
\node at (3.5, 3.5){$q_1$};
\node [left] at (3,0){\footnotesize$C_2(2,1)$};
\node at (3.5, -3.5){$p_2$};
\node [above] at (0,-3){\footnotesize$C_2(2,2)$};
\node at (-3.5, -3.5){$q_2$};
\node [right] at (-3,0){\footnotesize$C_2(1,2)$};
\node at (-3.5, 3.5){$p_1$};

\draw (-6,6)--(3, 3);
\node at (-2, 5.25){\footnotesize $C_2(0,1)$};
\draw (-6,6)--(-3, -3);
\node at (-4.25, 4.25){\footnotesize$C_2(0,2)$};
\draw (6,6)--(-3, 3);
\node at (2, 5.25){\footnotesize$C_2(1,0)$};
\draw (6,6)--(3, -3);
\node at (4.25, 4.25){\footnotesize$C_2(2,0)$};
\draw (6,-6)--(3, 3);
\node at (4.25, -4.24){\footnotesize$C_2(3,1)$};
\draw (6,-6)--(-3, -3);
\node at (2.25, -5.25){\footnotesize$C_2(3,2)$};
\draw (-6,-6)--(-3, 3);
\node at (-4.25, -4.25){\footnotesize$C_2(1,3)$};
\draw (-6,-6)--(3, -3);
\node at (-1.75, -5.25){\footnotesize$C_2(2,3)$};
\end{tikzpicture}

\vspace{0.15cm}

\begin{tikzpicture}[xscale=.5,yscale=.5]
\draw (-6,6)--(6,6)--(6,-6)--(-6,-6)--(-6,6);
\foreach \Point in {(-6,6),(6,6),(6,-6),(-6,-6)}{
    \node at \Point {\textbullet};
}
\node [above] at (0,6){\footnotesize$C_3(0,0)$};
\node at (-6.5,6.5){$r_0$};
\node [right] at (6,0){\footnotesize$C_3(3,0)$};
\node at (6.5,6.5){$s_0$};
\node [below] at (0,-6){\footnotesize$C_3(3,3)$};
\node at (6.5,-6.5){$r_3$};
\node [left] at (-6,0){\footnotesize$C_3(0,3)$};
\node at (-6.5,-6.5){$s_3$};

\draw (-3,3)--(3,3)--(3,-3)--(-3,-3)--(-3,3);
\foreach \Point in {(-3,3),(3,3),(3,-3),(-3,-3)}{
    \node at \Point {\textbullet};
}
\node [below] at (0,3){\footnotesize$C_3(1,1)$};
\node at (3.5, 3.5){$s_1$};
\node [left] at (3,0){\footnotesize$C_3(2,1)$};
\node at (3.5, -3.5){$r_2$};
\node [above] at (0,-3){\footnotesize$C_3(2,2)$};
\node at (-3.5, -3.5){$s_2$};
\node [right] at (-3,0){\footnotesize$C_3(1,2)$};
\node at (-3.5, 3.5){$r_1$};

\draw (-6,6)--(3, 3);
\node at (-2, 5.25){\footnotesize $C_3(0,1)$};
\draw (-6,6)--(-3, -3);
\node at (-4.25, 4.25){\footnotesize$C_3(0,2)$};
\draw (6,6)--(-3, 3);
\node at (2, 5.25){\footnotesize$C_3(1,0)$};
\draw (6,6)--(3, -3);
\node at (4.25, 4.25){\footnotesize$C_3(2,0)$};
\draw (6,-6)--(3, 3);
\node at (4.25, -4.25){\footnotesize$C_3(3,1)$};
\draw (6,-6)--(-3, -3);
\node at (2.25, -5.25){\footnotesize$C_3(3,2)$};
\draw (-6,-6)--(-3, 3);
\node at (-4.25, -4.25){\footnotesize$C_3(1,3)$};
\draw (-6,-6)--(3, -3);
\node at (-1.75, -5.25){\footnotesize$C_3(2,3)$};
\end{tikzpicture}

\caption{SLICE 1 of lines on $X_1$. These three diagrams depict intersection points within families. Intersections indicated with ``\textbullet''.} 
\label{fig:X_1-lines-slice1}
\end{figure}

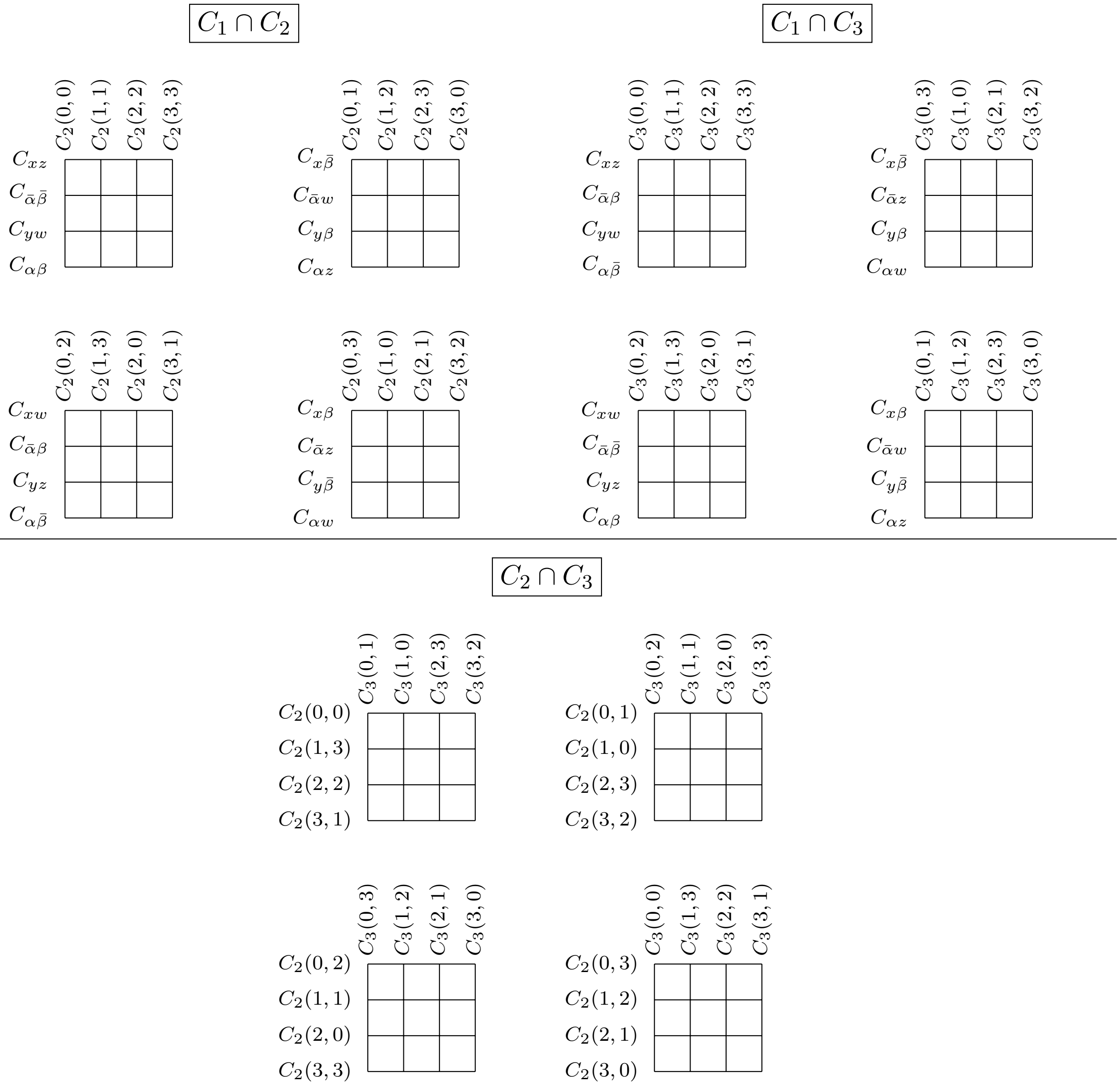
\begin{figure}
    \centering


\newcommand{\fourgrid}[4]{%
   \draw[step=1.0,black,thin] (0,0) grid (3,3);

  \foreach \a/\b [count=\k from 0] in {#3}{%
   
    \node[anchor=east,font=\scriptsize] at (-0.18,3-\k)
      {$#1(\a,\b)$};
  }

  \foreach \a/\b [count=\k from 0] in {#4}{%
      \node[
      rotate=90,
      font=\scriptsize
    ] at (\k,4.25) {$#2(\a,\b)$};
  }
}

\newcommand{\Cgrid}[4]{%
   \draw[step=1.0,black,thin] (0,0) grid (3,3);

  \foreach \a/\b [count=\k from 0] in {#3}{%
   
    \node[anchor=east,font=\scriptsize] at (-0.18,3-\k)
      {${#1}_{\a \b}$};
  }

  \foreach \a/\b [count=\k from 0] in {#4}{%
      \node[
      rotate=90,
      font=\scriptsize
    ] at (\k,4.25) {$#2(\a,\b)$};
  }
}

\begin{tikzpicture}[x=.5cm,y=.5cm]

\node[draw,inner sep=3pt] at (5,14.8) {$C_1\cap C_2$};

\begin{scope}[shift={(0,8)}]
\Cgrid{C}{C_2}
  {x/z,\bar{\alpha}/\bar{\beta},y/w,\alpha/\beta}
  {0/0,1/1,2/2,3/3}
\end{scope}

\begin{scope}[shift={(8,8)}]
 \Cgrid{C}{C_2}
   {x/\bar{\beta},\bar{\alpha}/w,y/\beta,\alpha/z}
   {0/1,1/2,2/3,3/0}
 \end{scope}

 \begin{scope}[shift={(0,1)}]
 \Cgrid{C}{C_2}
   {x/w,\bar{\alpha}/\beta,y/z,\alpha/\bar{\beta}}
   {0/2,1/3,2/0,3/1}
 \end{scope}

 \begin{scope}[shift={(8,1)}]
 \Cgrid{C}{C_2}
   {x/\beta,\bar{\alpha}/z,y/\bar{\beta},\alpha/w}
   {0/3,1/0,2/1,3/2}
 \end{scope}

 \node[draw,inner sep=3pt] at (21,14.8) {$C_1\cap C_3$};

 \begin{scope}[shift={(16,8)}]
 \Cgrid{C}{C_3}
   {x/z,\bar{\alpha}/\beta,y/w,\alpha/\bar{\beta}}
   {0/0,1/1,2/2,3/3}
 \end{scope}

 \begin{scope}[shift={(24,8)}]
 \Cgrid{C}{C_3}
   {x/\bar{\beta},\bar{\alpha}/z,y/\beta,\alpha/w}
   {0/3,1/0,2/1,3/2}
 \end{scope}

 \begin{scope}[shift={(16,1)}]
 \Cgrid{C}{C_3}
   {x/w,\bar{\alpha}/\bar{\beta},y/z,\alpha/\beta}
   {0/2,1/3,2/0,3/1}
 \end{scope}

 \begin{scope}[shift={(24,1)}]
 \Cgrid{C}{C_3}
   {x/\beta,\bar{\alpha}/w,y/\bar{\beta},\alpha/z}
   {0/1,1/2,2/3,3/0}
 \end{scope}
\end{tikzpicture}

\hrule
\vspace{.25cm}

\begin{tikzpicture}[x=.5cm,y=.5cm]

 \node[draw,inner sep=3pt] at (5,14.8) {$C_2\cap C_3$};

 \begin{scope}[shift={(0,8)}]
 \fourgrid{C_2}{C_3}
   {0/0,1/3,2/2,3/1}
   {0/1,1/0,2/3,3/2}
 \end{scope}

 \begin{scope}[shift={(8,8)}]
 \fourgrid{C_2}{C_3}
   {0/1,1/0,2/3,3/2}
   {0/2,1/1,2/0,3/3}
 \end{scope}

 \begin{scope}[shift={(0,1)}]
 \fourgrid{C_2}{C_3}
   {0/2,1/1,2/0,3/3}
   {0/3,1/2,2/1,3/0}
 \end{scope}

 \begin{scope}[shift={(8,1)}]
 \fourgrid{C_2}{C_3}
   {0/3,1/2,2/1,3/0}
   {0/0,1/3,2/2,3/1}
 \end{scope}

\end{tikzpicture}

\caption{SLICE 2 of lines on $X_1$ depicting intersections between families.}
\label{fig:X_1-lines-slice2}

\end{figure}

The bilinear form on the Picard lattice is determined by the intersection between these lines, which one can compute. To visualize these lines and their intersections, we display two ``slices'' of the full intersection diagram, displayed in Figures~\ref{fig:X_1-lines-slice1}, where we depict intersections within families and \ref{fig:X_1-lines-slice2}, where we depict intersections between families.  Some explanation of the figures is in order.

In Figure~\ref{fig:X_1-lines-slice2}, intersections on $X_1$ are simply depicted as the intersection of the lines in the figure. In Figure~\ref{fig:X_1-lines-slice1}, intersections are shown with ``\textbullet''. 
Where two lines in the figure intersect without ``\textbullet'', there is no intersection on $X_1$.

It is arranged so that the intersection points from the two slices are disjoint. That is to say, no intersection point between these lines joins lines from the same family and lines from a different family. Furthermore, all relevant lines which share the same intersection point are included in the same slice, and every intersection point between any two lines is included in one of the slices. Together, this ensures all information about intersections between relevant curves is faithfully represented.

Interaction between families depicted in Figure~\ref{fig:X_1-lines-slice2} forms a simple structure. Given two families, the lines in each family may be sorted into four groups, such that one group from the first family intersects exactly one group from the other family. Moreover, intersections between these groups are complete, in the sense that every line in the first group intersects every line in the second, with exactly two lines meeting at an intersection. From the diagrams, we can build the intersection matrix.  

We need to name the intersection points in Figure~\ref{fig:X_1-lines-slice1}. Let $t_x=[0:1:0:0]$ denote the intersection of the four curves $C_{x\ell}$ for $\ell\in\{z, w, \beta, \bar{\beta}\}$. Similarly, we denote the other intersection points between the $C_1$ family as $t_y, t_z, t_w, t_\alpha, t_{\bar\alpha}, t_\beta, t_{\bar\beta}$, in the same manner. 
Denote $p_m$ for the intersection between the four curves $C_2(m,\cdot)$, $q_n$ for the intersection of the four curves $C_2(\cdot, n)$, $r_m$ the intersection between the four curves $C_3(m,\cdot)$, $s_n$ for the intersection of the four curves $C_3(\cdot, n)$.


\subsubsection{Diagonal Automorphisms}

The group of automorphisms of $X_i$ induced by diagonal symmetries is
\[\overline{\Gdiag_i} := \Gdiag/G_i~,\]
with $G_i$ defined in \eqref{eq:definition-of-G2}. Note that any two diagonal symmetries that differ by an element of $G_i$ will produce the same automorphism of the the quotient by $G_i$, and each will extend to the minimal resolution.

Recall $\Gdiag$ is generated by order 8 elements $\sigma_L, \sigma_R$ defined in \eqref{eq:definition-sigmaL} and \eqref{eq:definition-sigmaR}. 
We let 
\[
\tau := \sigma_L + \sigma_R = \left(\frac{3}{8}, -\frac{1}{8}, \frac{3}{8}, -\frac{1}{8}\right).
\]
One can check that $2\tau\in J_W=G_1$, so in the quotient by $G_1$ we obtain the relation $2(\sigma_L +\sigma_R) = 2\tau = 0$. It turns out this gives us all of the relations. 

We can use this to give a standard set of generators for $\overline{G^{diag}_i}$. If we define
\[
A := \begin{pmatrix}
    8 & 0 & 2 \\
    0 & 8 & 2
\end{pmatrix}~,
\]
then
\[\overline{\Gdiag_1} = \Gdiag/G_1 \cong \ZZ^2/\image A\]
via $(m\sigma_L + n\sigma_R) + G_1 \mapsto (m, n) + \image A$. We have the following Smith decomposition,
\[
PAQ = 
\begin{pmatrix}
    1 & 0 \\
    -1 & 1
\end{pmatrix}
\begin{pmatrix}
    8 & 0 & 2 \\
    0 & 8 & 2
\end{pmatrix}
\begin{pmatrix}
    0 & 0 & -1 \\
    0 & 1 & -1 \\
    1 & 0 & 4
\end{pmatrix}
=
\begin{pmatrix}
    2 & 0 & 0 \\
    0 & 8 & 0
\end{pmatrix}
= D~,
\]
defining respectively the matrices $P, Q$ and $D$. The matrix $P$ represents an automorphism on $\ZZ^2$ and $Q$ an automorphism on $\ZZ^3$. Under the automorphism $P$, we have
\[
\image A \xrightarrow{P} \image(PA) = \image (PAQ) = \image D
\]
(using the fact that $Q$ is an automorphism). Hence, $P$ induces an isomorphism
\[
    \ZZ^2/\image A \xrightarrow{\overline{P}} \ZZ^2/\image D~.
\]

The standard generator $e_i + \image D$ for $\ZZ^2/\image D$ in terms of the generators on $\ZZ^2/\image A$ can be written $P^{-1}e_i + \image A$, or the $i$th column of $P^{-1}$. Since
\[
 P^{-1} =
 \begin{pmatrix}
    1 & 0 \\
    1 & 1 \\
 \end{pmatrix}~,
\]
then these generators are $(\sigma_L+\sigma_R) + G_1 = \tau + G_1$ and $\sigma_R + G_1$. We conclude
\[
\Gdiag/G_1 \cong \ZZ_2 \oplus \ZZ_8,
\]
with isomorphism given by 
\[
(m\tau + n\sigma_R) + G_1 \mapsto (m, n).
\]
This gives us a standard representation for the elements of $\overline{\Gdiag_1}$. 

Considering each of these symmetries as a diagonal matrix, and using the description of the volume form for a hypersurface, gives an easy way to check if a diagonal symmetry is symplectic. It is symplectic if the determinant of the corresponding diagonal matrix is 1, and it is nonsymplectic otherwise (see e.g. \cite[Proposition 2]{ABS}). For the sixteen elements of $\Gdiag/G_1$, the determinants are given in Table~\ref{table:determinant}, from which we can determine which are symplectic, and which are nonsymplectic.  

\begin{table}
    \renewcommand{\arraystretch}{1.25}
    \begin{tabular}{|c|c|c||c|c|c|}
        \hline
        $\bm{(m, n)}$ & Order & Determinant & $\bm{(m, n)}$ & Order & Determinant\\ \hline \hline
        (0, 0) & 1 & 1 & (1,0) & 2 & $-1$ \\ \hline
        (0, 1) & 8 & $i$ & (1,1) & 8 & $-i$ \\ \hline
        (0, 2) & 4 & $-1$ & (1,2) & 4 & $1$ \\ \hline
        (0, 3) & 8 & $-i$ & (1,3) & 8 & $i$ \\ \hline
        (0, 4) & 2 & $1$ & (1,4) & 2 & $-1$ \\ \hline
        (0, 5) & 8 & $i$ & (1,5) & 8 & $-i$ \\ \hline
        (0, 6) & 4 & $-1$ & (1,6) & 4 & $1$ \\ \hline
        (0, 7) & 8 & $i$ & (1,7) & 8 & $i$ \\ \hline
    \end{tabular}
    \vspace{0.1cm}
    \caption{The order and determinant of each element of $G^{diag}/G_1$.}
    \label{table:determinant}
\end{table}

We can compare the table with the 
exact sequence \eqref{eq:automorphism-short-exact-sequence}. In our example, this becomes
\[
0\to \ZZ/4\ZZ \to \ZZ/2\ZZ \oplus \ZZ/8\ZZ \to \ZZ/4\ZZ \to 0.
\]
So we see there are four symplectic automorphisms, and and any power of any automorphism higher than four will yield a symplectic automorphism. 

\subsection{Action of $\sigma_R$ and $\sigma_L$ on $X_1$}
In order to make other computations much simpler, we next describe the action of $\sigma_R$ and $\sigma_L$ on the relevant curves, which the interested reader can easily verify.  

One should notice the following: Suppose we have a linear form $ax+by+cz+dw$. If $[x_0:y_0:z_0:w_0]\in V(ax+by+cz+dw)$, and we have a diagonal symmetry (written multiplicatively) $(g_x,g_y, g_z, g_w)$ acting, then $[g_xx_0:g_yy_0:g_zz_0:g_ww_0]\in V((ag_x^{-1})x+(bg_y^{-1})y+(cg_z^{-1})z+(dg_w^{-1})w)$. So for example $\sigma_L$ will send $C_{\alpha\beta}$ to $V(\zeta^{-3}x+i\zeta y, \beta)=V(x-iy, \beta)=C_{\bar{\alpha}\beta}$, where as before $\zeta$ is a primitive 8th root of unity. 

With this in mind, we can describe the action of $\sigma_R$ and $\sigma_L$. First, we will describe the action of $\sigma_R$. Note that $\sigma_R$ preserves the zero loci of the linear forms $x$, $y$, $z$, $w$, $\alpha$ and $\bar{\alpha}$, but it swaps $V(\beta)$ and $V(\bar{\beta})$. Thus, it leaves $C_{xz}$, $C_{yz}$, $C_{xw}$, $C_{yw}$, $C_{\alpha z}$, $C_{\alpha w}$, $C_{\bar\alpha z}$, and $C_{\bar\alpha w}$ invariant, but it exchanges
\begin{align*}
C_{x\beta} &\leftrightarrow C_{x\bar\beta}  \quad &C_{\alpha\beta} \leftrightarrow C_{\alpha\bar\beta}\\
C_{y\beta} &\leftrightarrow C_{y\bar\beta} \quad &C_{\bar\alpha\beta} \leftrightarrow C_{\bar\alpha\bar\beta}.
\end{align*}
Furthermore, it fixes $t_x,t_y, t_z, t_w, t_\alpha$, and $t_{\bar{\alpha}}$, but it exchanges $t_{\beta}$ and $t_{\bar\beta}$. 

For the other curves and intersection points, we have:
\begin{align}
    &C_2(m,n)\to C_3(m+3,n+3)\nonumber\\
    &C_3(m,n)\to C_2(m+2,n+2)\nonumber\\
    &p_m\to r_{m+3}\label{eq:sigmaR}\\
    &q_n\to s_{n+3}\nonumber\\
    &r_m\to p_{m+2}\nonumber\\
    &s_n\to q_{n+2}~.\nonumber
\end{align}

Next we describe the action of $\sigma_L$. It preserves the zero locus of $x$, $y$, $z$, $w$, $\beta$ and $\bar{\beta}$, but it swaps $V(\alpha)$ and $V(\bar{\alpha})$. Thus, it leaves 
$C_{xz}$, $C_{yz}$, $C_{xw}$, $C_{yw}$, $C_{x\beta}$, $C_{x\bar\beta}$, $C_{y\beta}$, and $C_{y\bar\beta}$ invariant, but it performs the following exchanges
\begin{align*}
C_{\alpha z} &\leftrightarrow C_{\bar\alpha z}  \quad &C_{\alpha\beta} \leftrightarrow C_{\bar\alpha\beta}\\
C_{\alpha w} &\leftrightarrow C_{\bar\alpha w} \quad &C_{\alpha\bar\beta} \leftrightarrow C_{\bar\alpha\bar\beta}.
\end{align*}
Furthermore, $\sigma_L$ fixes $t_x,t_y, t_z, t_w, t_{\beta}$ and $t_{\bar\beta}$, but it exchanges $t_\alpha$ and $t_{\bar{\alpha}}$. 

For the other curves and intersection points, $\sigma_L$ acts as follows:
\begin{align}
    &C_2(m,n)\to C_3(n+2,m+2)\nonumber\\
    &C_3(m,n)\to C_2(n+1,m+1)\nonumber\\
    &p_m\to s_{m+2}\label{eq:sigmaL}\\
    &q_n\to r_{n+2}\nonumber\\
    &r_m\to q_{m+1}\nonumber\\
    &s_n\to p_{n+1}~.\nonumber
\end{align}

Since $\sigma_R$ and $\sigma_L$ generate every diagonal automorphism, we can easily track the action of any diagonal automorphism by using the actions in \eqref{eq:sigmaR} and \eqref{eq:sigmaL} given above.

\subsection{Surface $X_4$}

Recall that $X_4$ is the BHK mirror of $X_1$. We now turn our attention to this surface.

\subsubsection{Geometry}\label{ss:X4geom}

Recall that $\overline{G}_4=G_4/J_W$ is generated by $\sigma_L-\sigma_R$. We will now describe the action of $\sigma_L-\sigma_R$ on the lines on $X_1$, to help us find a generating set for $\Pic(X_4)$. 
First notice that $\sigma_L-\sigma_R$ fixes the zero locus of the linear forms $x,y,z,w$, but it exchanges the zero locus of $\alpha$ with $\bar\alpha$ and the zero locus of $\beta$ with $\bar\beta$. Thus, $\sigma_L-\sigma_R$ keeps $C_{xz}$, $C_{yz}$, $C_{xw}$, and $C_{yw}$ invariant (though not fixed pointwisely), and exchanges the following: 
\begin{align*}
C_{x\beta} &\leftrightarrow C_{x\bar\beta}  \quad &C_{\alpha\beta} \leftrightarrow C_{\bar\alpha\bar\beta}\\
C_{y\beta} &\leftrightarrow C_{y\bar\beta} \quad &C_{\bar\alpha\beta} \leftrightarrow C_{\alpha\bar\beta}.\\
C_{\alpha z} &\leftrightarrow C_{\bar\alpha z}  \quad &\\
C_{\alpha w} &\leftrightarrow C_{\bar\alpha w} \quad &
\end{align*}

Referring back to Figure \ref{fig:X_1-lines-slice1}, we observe the points labeled $t_x, t_y, t_z, t_w$ are fixed, and there are two equivalence classes among the other four points $t_\alpha\sim t_{\bar{\alpha}}$ and $t_\beta\sim t_{\bar{\beta}}$. 
Therefore, in the quotient, there are four $A_3$ singularities (at the points labeled $t_x, t_y, t_z, t_w$) and two $A_1$ singularities (at the glued points $t_\alpha\sim t_{\bar{\alpha}}$ and $t_\beta\sim t_{\bar{\beta}}$). 
Since $\overline{G_4} \cong \ZZ/4\ZZ$, by \cite{Xiao} these are all the singularities resulting from this group action. Notice that the families $C_2(m,n)$ and $C_3(m,n)$ do not intersect these singular points, and thus are not effected by the resulting blowup.

Next consider the action of $\sigma_L-\sigma_R$ on the other two families of lines. Using the action described in \eqref{eq:sigmaR} and \eqref{eq:sigmaL}
\begin{align*}
    C_2(m,n) &\mapsto C_2(n-1, m-1) \\
    C_3(m,n) &\mapsto C_3(n-1, m-1).
\end{align*}
We get the following eight equivalence classes under this action:
\begin{align*}
    [C_2(0,0)] &:= \{C_2(0,0),C_2(3,3), C_2(2,2), C_2(1,1)\} \\
    [C_2(0,2)] &:= \{C_2(0,2),C_2(3,1), C_2(2,0), C_2(1,3)\} \\
    [C_2(0,1)] &:= \{C_2(0,1),C_2(2,1), C_2(2,3), C_2(0,3)\}\\
    [C_2(1,0)] &:= \{C_2(1,0),C_2(1,2), C_2(3,2), C_2(3,0)\}\\
    [C_3(0,0)] &:= \{C_3(0,0),C_3(3,3), C_3(2,2), C_3(1,1)\} \\
    [C_3(0,2)] &:= \{C_3(0,2),C_3(3,1), C_3(2,0), C_3(1,3)\} \\
    [C_3(0,1)] &:= \{C_3(0,1),C_3(2,1), C_3(2,3), C_3(0,3)\}\\
    [C_3(1,0)] &:= \{C_3(1,0),C_3(1,2), C_3(3,2), C_3(3,0)\},
\end{align*}
as well as the following equivelence classes of intersection points:
\begin{align*}
\{p_0, q_1,p_2,q_3\},\quad\{p_1,q_2,p_3,q_0\},\quad\{r_0, s_1,r_2,s_3\},\quad\{r_1,s_2,r_3,s_0\}.
\end{align*}
This leaves us with 8 lines and four intersection points between these lines. 

\begin{figure}[ht]
\centering\begin{tikzpicture}[xscale=.7,yscale=.5]
\foreach \Point in {(2, 4), (4, -4), (-4, 4), (-2, -4)}{
    \node at \Point {\small$\blacksquare$};
}
\foreach \Point in {(2, 2), (4, 2), (2, -2), (4, -2), (-2, 2), (-4, 2), (-2, -2), (-4, -2)}{
    \node at \Point {\textbullet};
}

\draw (-4,2)--(-2,2);
\node [left] at (-4,2){$C_{xz}$};
\draw (-4,-2)--(-2,-2);
\node [left] at (-4,-2){$C_{yw}$};
\draw (2,2)--(4,2);
\node [right] at (4,2){$C_{xw}$};
\draw (2, -2)--(4,-2);
\node [right] at (4,-2){$C_{yz}$};

\draw (-4,4)--(2,4);
\node [above] at (-1,4){$C_{\alpha\bar{\beta}\sim \bar{\alpha}\beta}$};
\draw (-2,-4)--(4,-4);
\node [below] at (1,-4){$C_{\alpha\beta\sim \bar{\alpha}\bar{\beta}}$};

\draw (-4,4)--(-4,-2);
\node [left] at (-4, 0){\small$[C_3(0,0)]$};
\node [above] at (-4, 4){$\times 2$};
\draw (-2,2)--(-2,-4);
\node [below] at (-2, -4){$\times 2$};
\node [right] at (-2, -.5){\small$[C_2(0,0)]$};

\draw (2,4)--(2,-2);
\node [above] at (2, 4){$\times 2$};
\node [left] at (2, .5){\small$[C_2(0,2)]$};
\draw (4,2)--(4,-4);
\node [below] at (4, -4){$\times 2$};
\node [right] at (4, 0){\small$[C_3(0,2)]$};

\end{tikzpicture}
\caption{SLICE 2 of lines on $X_4$. Intersections indicated with ``\textbullet''. Two distinct intersections between the same lines indicated with $\times 2$.}
\label{fig:X_4-lines-slice2}
\end{figure}
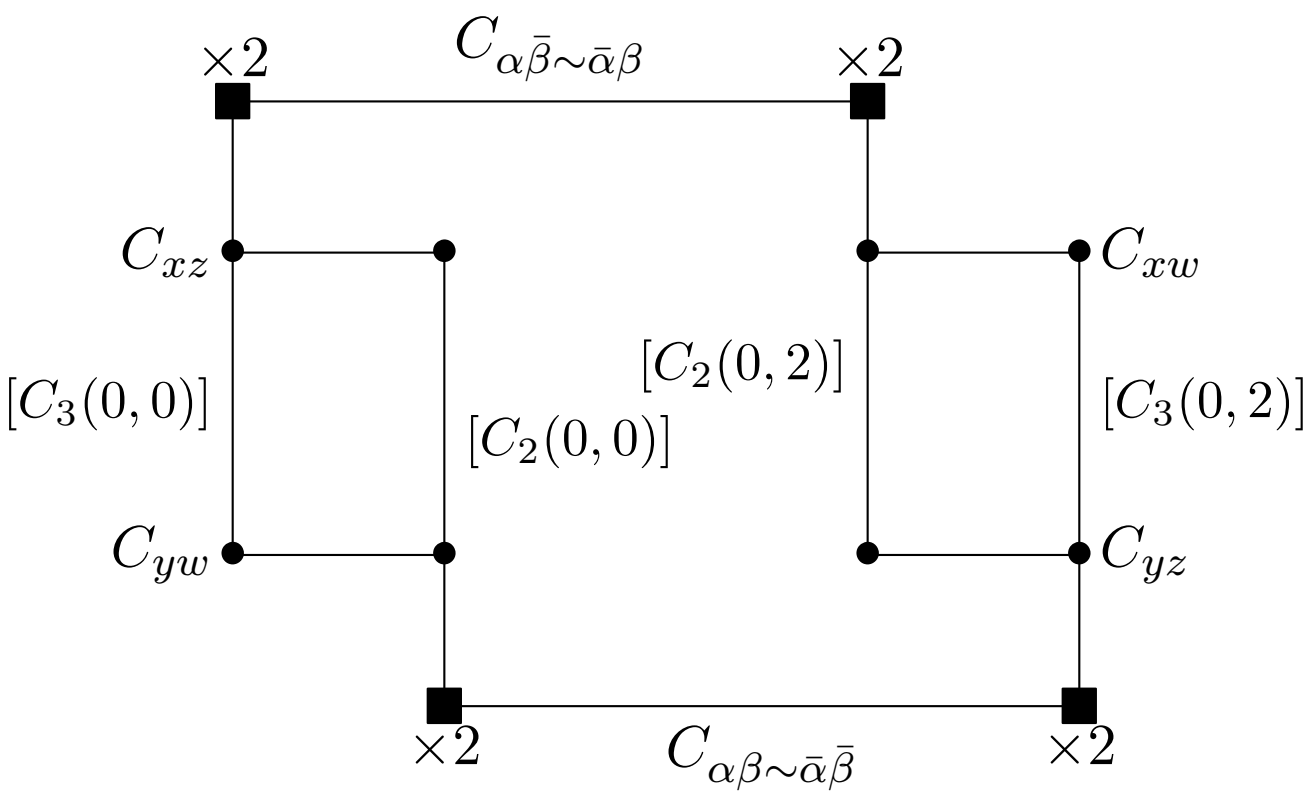

Resolving the singularities, we obtain the K3 surface $X_4$. 
We can represent the collection of curves in Figures~\ref{fig:X_4-lines-slice1} and \ref{fig:X_4-lines-slice2}, including the exceptional curves $E_i$ resulting from the minimal resolution of the four $A_3$ singularities and two $A_1$ singularities.\footnote{In order to see that the blowup of each $A_3$ singularity results in the depicted configuration of curves, one can first construct $X_2$, and then act by a $\ZZ/2\ZZ$ group action on the resulting surface to construct $X_4$. We have chosen to present $X_4$ first for clarity of exposition.} Furthermore, one can check that the remaining intersections between the relevant curves is as depicted in Figure~\ref{fig:X_4-lines-slice2}. Using the same convention as before, a ``\textbullet'' denotes an intersection point in both figures. We have used a square node labeled ``$\times 2$'' to indicate that certain curves intersect in two distinct points. The image of the curves $C_2(0,1)$, $C_2(1,0)$, $C_3(1,0)$, and $C_3(0,1)$ become nodal curves, so they are not depicted in either diagram, and we will not consider them further.

In Section~\ref{sss:X4invariantLattices}, we will show that these curves generated the Picard lattice $\Pic X_4$.

\begin{figure}[ht]
\centering\begin{tikzpicture}[xscale=1.25,yscale=1.25]
\foreach \Point in {(1.,2.41421),(2.,3.41421),(1.,-2.41421),(2.,-3.41421),(-1.,2.41421),(-2.,3.41421),(-1.,-2.41421),(-2.,-3.41421),(2.41421,1.),(3.41421,2.),(2.41421,-1.),(3.41421,-2.),(-2.41421,1.),(-3.41421,2.),(-2.41421,-1.),(-3.41421,-2.),(2.70711,2.70711),(2.70711,-2.70711),(-2.70711,2.70711),(-2.70711,-2.70711),(1.,1.),(1.,-1.),(-1.,1.),(-1.,-1.),(1,0.333333),(1,-0.333333),(-1,0.333333),(-1,-0.333333)}{
    \node at \Point {\textbullet};
}

\node at (1.3,3.11421){$E_{xz}$};
\draw (1.,2.41421)--(2.,3.41421);
\node at (3.11421,1.3){$E_{xw}$};
\draw (2.41421,1.)--(3.41421,2.);
\node at (2.90711,2.90711){$E_{x}$};
\draw[dashed] (2.,3.41421)--(3.41421,2.);
\node at (2.3,1.70711){$C_{x,\beta \sim \bar{\beta}}$};
\draw (1, 1)--(2.70711,2.70711);
\node [left] at (2.41421,-.4){$C_{xw}$};
\draw[dashed] (2.41421,1.)--(2.41421,-1.);

\node at (3.11421,-1.3){$E_{wx}$};
\draw (2.41421,-1.)--(3.41421,-2.);
\node at (1.3,-3.11421){$E_{wy}$};
\draw (1.,-2.41421)--(2.,-3.41421);
\node at (2.90711,-2.90711){$E_{w}$};
\draw[dashed] (3.41421,-2.)--(2.,-3.41421);
\node at (2,-2){$C_{w,\alpha \sim \bar{\alpha}}$};
\draw (-1, -1)--(2.70711,-2.70711);
\node [below] at (0,-2.41421){$C_{yw}$};
\draw[dashed] (1,-2.41421)--(-1,-2.41421);

\node at (-1.3,-3.11421){$E_{yw}$};
\draw (-1.,-2.41421)--(-2.,-3.41421);
\node at (-3.11421,-1.3){$E_{yz}$};
\draw (-2.41421,-1.)--(-3.41421,-2.);
\node at (-2.90711,-2.90711){$E_{y}$};
\draw[dashed] (-2.,-3.41421)--(-3.41421,-2.);
\node at (-2,-2){$C_{y,\beta \sim \bar{\beta}}$};
\draw (1, -1)--(-2.70711,-2.70711);
\node [right] at (-2.41421,-.4){$C_{yz}$};
\draw[dashed] (-2.41421,-1.)--(-2.41421,1.);

\node at (-3.11421,1.3){$E_{zy}$};
\draw (-2.41421,1.)--(-3.41421,2.);
\node at (-1.3,3.11421){$E_{zx}$};
\draw (-1.,2.41421)--(-2.,3.41421);
\node at (-2.90711,2.90711){$E_{z}$};
\draw[dashed] (-3.41421,2.)--(-2.,3.41421);
\node at (-2.3,1.70711){$C_{z,\alpha \sim \bar{\alpha}}$};
\draw (-1, 1)--(-2.70711,2.70711);
\node [above] at (0,2.41421){$C_{xz}$};
\draw[dashed] (-1,2.41421)--(1,2.41421);

\node [left] at (-1,0){$E_{\alpha \sim \bar{\alpha}}$};
\draw[dashed] (-1,-1)--(-1,1);
\node [right] at (1,0){$E_{\beta \sim \bar{\beta}}$};
\draw[dashed] (1,-1)--(1,1);
\node [above] at (0, 0.333333){$C_{\alpha\beta \sim \bar{\alpha}\bar{\beta}}$};
\draw (-1,0.333333)--(1,0.333333);
\node [below] at (0, -0.333333){$C_{\alpha\bar{\beta} \sim \bar{\alpha}\beta}$};
\draw (-1,-0.333333)--(1,-0.333333);

\foreach \Point in {(4.5, 0.5), (4.5, -0.5)}{
    \node at \Point {\textbullet};
}
\draw (4.5, -0.5) to[bend left = 45] node[left]{\small$[C_2(0,0)$]} (4.5,0.5);
\draw (4.5,-0.5) to[bend right=45] node[right] {\small$[C_2(0,2)]$} (4.5,0.5);

\foreach \Point in {(-4.5, 0.5), (-4.5, -0.5)}{
    \node at \Point {\textbullet};
}
\draw (-4.5,-0.5) to[bend left=45]  node[left]  {\small$[C_3(0,2)]$} (-4.5,0.5);
\draw (-4.5,-0.5) to[bend right=45] node[right] {\small$[C_3(0,0)]$} (-4.5,0.5);

\end{tikzpicture}
\caption{SLICE 1 of lines on $X_4$. Intersections indicated with ``\textbullet''. The dashed lines represent lines that are fixed under the action of the nonsymplectic involution $2\sigma_R$.}
\label{fig:X_4-lines-slice1}
\end{figure}
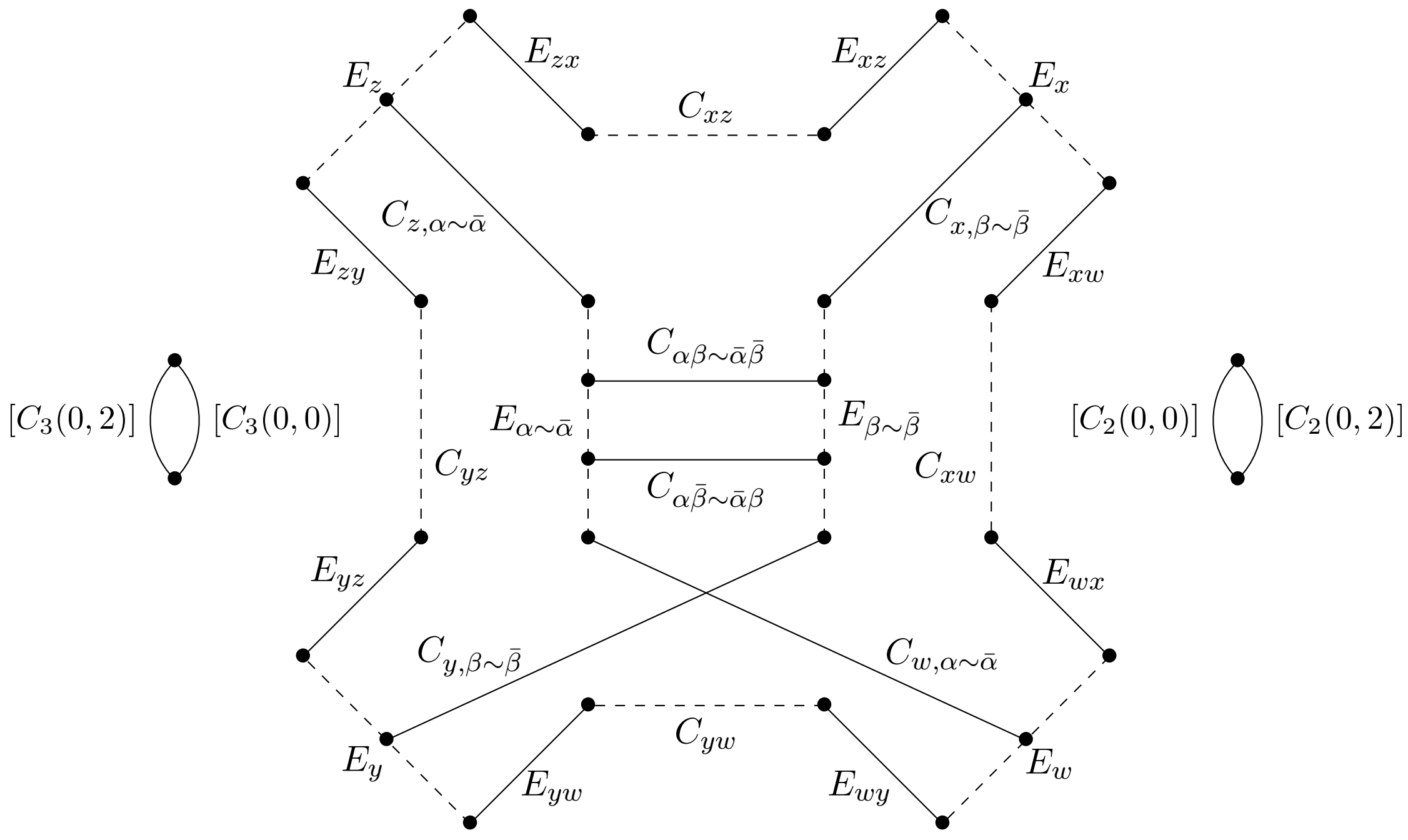

\subsubsection{Diagonal Automorphisms}

The diagonal automorphisms are
\[
\overline{\Gdiag} := \Gdiag / G_4~,\textrm{ with}\quad G_4 =\left\langle j_W, \sigma_L-\sigma_R\right\rangle~.\]
We have $\sigma_L + G_4 = \sigma_R + G_4$, so any element of $\Gdiag / G_4$ may be written as $n\sigma_R + G_4$ for some integer $n$. Moreover, one may easily show $\sigma_R + G_4$ is order 4, so
\[
\overline{\Gdiag_4} \cong \ZZ_4
\]
given by 
\[
n\sigma_R + G_4 \mapsto n~.
\]
In fact, from Table~\ref{table:determinant}, we see that every element is nonsymplectic and since $\det \sigma_R=i$, in fact $\sigma_R$ is purely nonsymplectic of order 4.

\subsection{Invariant Sublattices of the $X_1$, $X_4$ Pair}
As we have seen previously, $X_1$ and $X_4$ are a BHK mirror pair. According to Conjecture~\ref{conj:main}, we should be able to polarize $X_1$ by the invariant lattice of some nonsymplectic diagonal automorphism $\sigma\in \overline{\Gdiag_1}$ and $X_4$ by another nonsymplectic diagonal automorphism $\sigma^\vee\in \overline{\Gdiag_4}$ so that $S(\sigma^\vee)\cong S(\sigma)^\vee$. 

As in equation~\eqref{eq:mirrorRank}, a necessary condition for these two lattices to be mirror lattices is
\[
\rk S(\sigma^\vee)=20-\rk S(\sigma).
\]
Thus the next step is to compute the rank of the invariant lattice for each nonsymplectic element of $\overline{\Gdiag_1}$ and $\overline{\Gdiag_4}$, respectively, and then look for two lattices that satisfy the rank condition above. We will first make two important reductions. 

The first reduction is that we only need to compute the invariant lattice of a small subset of each group, that is, we only need to compute the invariant lattice of a single generator for each cyclic subgroup of $\overline{\Gdiag_i}$. Indeed, we first notice that if $\sigma$ is an automorphism, then $S(\sigma)\subseteq S(\sigma^k)$ for any $k\geq 1$. Furthermore, if $\sigma$ has order $n$ and $\gcd(n,k)=1$, then $\sigma^k$ generates the same cyclic subgroup as $\sigma$. Thus, $(\sigma^k)^m=\sigma$ for some $m\geq 1$, and hence $S(\sigma^k)\subseteq S(\sigma)$. Also, since the invariant lattice of a symplectic automorphism is not a sublattice of the Picard lattice, we only consider nonsymplectic diagonal symmetries. 

The second reduction is that it is enough for our purposes to compute the rank of each invariant lattice, rather than actually trying to identify the lattice. The following lemma gives us a simple way to check the rank. 

\begin{lemma}\label{lem:rankinvariant}
Let $G$ be a finite group and $L$ be a free $\ZZ$-module generated by a set $S\subseteq L$. Suppose that $G$ acts on $S$ (i.e. $S$ is $G$-stable), and let $M\subseteq L$ be the sublattice generated by the orbit sums of the action of $G$ on $S$. Then the invariant lattice of $G$ in $L$ has the same rank as $M$, i.e. $\rk M=\rk L^G$. 
\end{lemma}

In other words, we can compute the rank of the invariant lattice simply by computing the rank of the lattice generated by the orbit sums of the generators of $\Pic(X_i)$ for each diagonal symmetry. For $X_1$ our generating set is the set of 48 lines, and for $X_4$ this generating set consists of the strict transform of the orbits of lines, and the exceptional curves from the blow-up. This computation can be done quite quickly with a computer. 

Lemma~\ref{lem:rankinvariant} is a standard application of \cite[Prop. 3.3.1]{lorenz}. We will briefly outline the proof. 

\begin{proof}
Let $\ZZ[S]$ denote the free group generated by $S$. Then we have a surjective function $\ZZ[S]\to L$. After tensoring with $\QQ$, this extends to a $G$-equivariant surjective linear transformation 
\[
\phi: \QQ[S]\to L_\QQ:=L\otimes_\ZZ \QQ.
\]
Here $\QQ[S]$ denotes the vector space with $S$ as a basis. Now suppose $w\in (L_\QQ)^G$, and choose $v\in \QQ[S]$ such that $\phi(v)=w$. Notice also that $v_G:=\frac{1}{|G|}\sum_{g\in G} g(v)\in (\QQ[S])^G$ and that $\phi(v_G)=w$ (since $\phi$ is equivariant). 
Thus, we get a surjective linear transformation 
\[
(\QQ[S])^G\to (L_\QQ)^G.
\]
Since $\QQ[S]$ is a permutation representation, it is a well-known fact that $(\QQ[S])^G$ is generated by orbit sums. Hence, the orbit sums also generate $(L_\QQ)^G$. But the image of the orbit sums is just $M_\QQ$, i.e. $M_\QQ=(L_\QQ)^G$. By \cite[Prop. 3.3.1]{lorenz} we have $(L_\QQ)^G\cong L^G\otimes_\ZZ\QQ$. 
\end{proof}

\subsubsection{Diagonal Generator Action} \label{sec:diagonal-generator-action}

Notice that the $\tau$ and $\sigma_R$ cosets generate both $\overline{\Gdiag_1}$ and $\overline{\Gdiag_4}$ (though $\tau$ is trivial in the latter case). 

One can easily check that $\tau$ preserves the zero locus of the linear forms $x,y,z,w$, but permutes the zero locus of $\alpha$ and $\bar\alpha$ and $\beta$ and $\bar\beta$, resp., whereas $\sigma_R$ preserves all of these except for the zero loci of $\beta$ and $\bar\beta$, which are permuted. 

We conclude that for $\tau$, we get the following action on the $C_1$ family of lines: $C_{xz}, C_{xw}, C_{yz}, C_{yw}$ are $\tau$-invariant, and
\begin{align*}
    C_{x\beta} &\xleftrightarrow{\tau} C_{x\bar{\beta}} & C_{y\beta} &\xleftrightarrow{\tau} C_{y\bar{\beta}} & C_{\alpha\beta} &\xleftrightarrow{\tau}C_{\bar{\alpha}\bar{\beta}} \\
     C_{z\alpha} &\xleftrightarrow{\tau} C_{z\bar{\alpha}} & C_{w\alpha} &\xleftrightarrow{\tau} C_{w\bar{\alpha}} &  C_{\alpha\bar{\beta}} &\xleftrightarrow{\tau}C_{\bar{\alpha}\beta}~.
\end{align*}
Furthermore, we have the following action on the families $C_2$ and $C_3$, 
\begin{align*}
    C_2(m,n) &\mapsto C_2(n, m)\\
    C_3(m,n) &\mapsto C_3(n, m).
\end{align*}

In the case of $\sigma_R$, we find $C_{xz}, C_{xw}, C_{yz}, C_{yw}, C_{\alpha z}, C_{\bar{\alpha}z}, C_{\alpha w}, C_{\bar{\alpha}w}$ are $\sigma_R$-invariant, and the following are exchanged:  
\begin{align*}
    C_{x\beta} &\xleftrightarrow{\sigma_R} C_{x\bar{\beta}} & C_{\alpha\beta} &\xleftrightarrow{\sigma_R}C_{\alpha\bar{\beta}} \\
    C_{y\beta} &\xleftrightarrow{\sigma_R} C_{y\bar{\beta}} & C_{\bar{\alpha}\beta} &\xleftrightarrow{\sigma_R}C_{\bar{\alpha}\bar{\beta}}~.
\end{align*}
Furthermore, we have the following action on the other two families of lines: 
\begin{align*}
    C_2(m,n) &\mapsto C_3(m+3, n+3)\\
    C_3(m,n) &\mapsto C_2(m+2, n+2).
\end{align*}

From these computations, it will be straightforward to describe the action of every nonsymplectic diagonal symmetry, which we do below. 

\subsubsection{Diagonal Invariant Sublattices on $X_1$}

From the previous section, we can quickly compute the orbits of lines under the action of each nonsymplectic automorphism in $\overline{\Gdiag_1}$. Recall, we can represent each automorphism as an element of $\ZZ/2\ZZ\times \ZZ/8\ZZ$ via $m\tau + n\sigma_R\mapsto (m,n)$. We only need to consider the following automorphisms (cf. Table~\ref{table:determinant}), as each cyclic subgroup is generated by one of these: 
\begin{itemize}
    \item First, we have $(0,1)$ and $(1,1)$, both are automorphisms of order 8, whose fourth power is symplectic; 
    \item Second, we have $(0,2)$ which has order 4, but its square is symplectic; 
    \item Finally we have two nonsymplectic involutions $(1,0)$ and $(1,4)$. 
\end{itemize}

The orbits of each automorphism acting on the set of lines on $X_1$ is given in the following table, together with the rank of the invariant lattice. We have adopted the notation that $\{\alpha,\bar\alpha\}$ represents all of the orbits of the form $\{C_{\alpha \ell}, C_{\bar\alpha \ell}\}$, and similarly for $\{\beta,\bar\beta\}$. For the automorphisms $(1,0)$ and $(1,4)$, since both $\{\alpha,\bar\alpha\}$ and $\{\beta,\bar\beta\}$ are listed, this includes the two orbits $\{C_{\alpha\beta},C_{\bar\alpha\bar\beta}\}$ and $\{C_{\bar\alpha\beta},C_{\alpha\bar\beta}\}$.

\begin{center}
    \renewcommand{\arraystretch}{1.25}
    \begin{tabular}{|c|l|c|}
        \hline
        $\bm{(m, n)}$ &  \textbf{Nontrivial Orbits} & \textbf{Sublattice Rank}\\ \hline \hline
        (0, 1) &  $\{\beta, \bar{\beta}\}, \{C_k(m+i,n+i)\}_{i\in\ZZ/4\ZZ,~k \in \{2, 3\}}$ & 7 \\ \hline
        (0, 2) &  $\{C_k(m+i,n+i)\}_{i\in\ZZ/4\ZZ}$ for $k = 2, 3$ & 10 \\ \hline
        (1, 0) &  $\{\alpha, \bar{\alpha}\}, \{\beta, \bar{\beta}\}, \{C_k(m,n),C_k(n,m)\}$ for $m \ne n$, $k = 2, 3$ & 12 \\ \hline
        (1, 1) &  $\{\alpha, \bar{\alpha}\}, \{C_2(m+i,n+i), C_3(n+i,m+i)\}_{i\in\ZZ/4\ZZ}$ & 7 \\ \hline
        (1, 4) &  $\{\alpha, \bar{\alpha}\}, \{\beta, \bar{\beta}\}, \{C_k(m,n),C_k(n+2,m)\}$ for $k = 2, 3$& 12 \\ \hline
    \end{tabular}
\end{center}

\subsubsection{Diagonal Invariant Sublattices on $X_4$}\label{sss:X4invariantLattices}

There are only two diagonal nonsymplectic automorphisms on $X_4$ that we need to consider, namely $\sigma_R$ and $2\sigma_R$. 
We would like to know that the set of lines described in Section~\ref{ss:X4geom} generates all of $\Pic(X_4)$. We can do this by considering the (purely) nonsymplectic involution $2\sigma_R=(\tfrac{1}{4}, \tfrac{1}{4},0,0)=(0,0,\tfrac{3}{4},\tfrac{3}{4})$.

In Section~\ref{sec:background}, we learn that the invariant lattice of any nonsymplectic involution can be determined by examining the topology of the fixed locus, and that this lattice is a primitive sublattice of $\Pic(X_4)$. One can check that $2\sigma_R$ fixes the dashed curves in Figure~\ref{fig:X_4-lines-slice1}. For example, it keeps $C_{xz}$, $C_{yz}$, and $C_{\alpha z\sim \bar\alpha z}$ invariant. Hence prior to the blowup, it fixes $t_z$. Thus, $2\sigma_R$ leaves the exceptional curves resulting from the blowup invariant, and therefore, $E_z$ must have at least three fixed points. As $E_z$ is a rational curve, it must be fixed. Similar reasoning will yield all of the fixed curves (see also \cite[Lemma 7.1]{dillies16}). 

So the fixed locus of $2\sigma_R$ consists of 10 rational curves. Comparing with Figure~\ref{fig:2elementary}, we can see that the invariant lattice $S(2\sigma_R)$ has rank 20 with discriminant group $(\ZZ/2\ZZ)^2$ and is primitively embedded in $\Pic X_4$. 
As $S(2\sigma_R)$ is primitively embedded in $\Pic(X_4)$ and both have rank 20, we have $\Pic(X_4)=S(2\sigma_R)$. 

Computationally, one can also check that the lattice $L$ generated by the list of curves described in Section~\ref{ss:X4geom} has rank 20. Furthermore, one can verify that $L$ also has discriminant group $(\ZZ/2\ZZ)^2$. Hence $L=\Pic(X_4)$.

Recall that in the quotient $\sigma_R\equiv\sigma_L$. We have already seen that $2\sigma_R$ fixes the entire Picard lattice of $X_4$, which means $\sigma_R$ will have orbits of size at most 2. As for the action of $\sigma_R$, we can see that each exceptional curve will be invariant, as well as all of the images of the original lines from family $C_1$, except for $C_{\alpha\beta\sim\bar\alpha\bar\beta}$ will be permuted with $C_{\bar\alpha\beta\sim\alpha\bar\beta}$. 

We also see the following families of curves identified under the action: 
\begin{align*}
    [C_2(0,0)] &\leftrightarrow [C_3(0,0)]\\
    [C_2(0,2)] &\leftrightarrow [C_3(0,2)]
\end{align*}

In the following table, we list the orbits of the action of $\sigma_R$ and the ranks of the respective invariant lattices. Here $(m,n)$ can only take on two values, either $(0,0)$ or $(0,2)$, as just described. 
\begin{center}
    \renewcommand{\arraystretch}{1.25}
    \begin{tabular}{|c|l|c|}
        \hline
        $\bm{n}$ &  \textbf{Nontrivial Orbits} & \textbf{Sublattice Rank}\\ \hline \hline
        1 &  $\{(\alpha\beta\sim\bar{\alpha}\bar{\beta}), (\alpha\bar{\beta} \sim \bar{\alpha}\beta)\}, \{C_2(m,n),C_3(m,n)\}$ & 19 \\ \hline
        2 &  {--} & 20 \\ \hline
    \end{tabular}
\end{center}
\bigskip

Based on Equation~\eqref{eq:mirrorRank}, we can see that no invariant lattice of a diagonal symmetry from $X_1$ will be mirror to either of these invariant lattices. Thus Conjecture 1 is not true for diagonal symmetries. But in the next section, we will consider two nondiagonal symmetries. 

\subsubsection{Non-Diagonal Invariant Sublattices}

From Section~\ref{ss:nondiagonal}, we can see that there are also certain nondiagonal symmetries of our K3 surfaces, coming from permuting the variables. Four of these are nonsymplectic, namely $(xzyw)$, $(xwyz)$, $(xy)$ and $(zw)$. 

The first two are inverse to each other, and so generate the same subgroup of $G^{max}$. Thus they have the same invariant lattice. Each of these is nonsymplectic of order 4, but their square is symplectic. 

The latter two are nonsymplectic involutions. However, there is a symmetry of $W$, from which we see that the invariant lattice for $(xy)$ will be isomorphic to the invariant lattice for $(zw)$. So we will only consider $(xy)$. 

We can compute the ranks of the invariant lattice via the orbit sums. These are given in the following tables. As before, we represents all of the orbits of the form $\{C_{\alpha \ell}, C_{\bar\alpha \ell}\}$, and similarly with $\{x,y\}$. This includes on $X_4$ orbits containing exceptional curves, such as $\{E_{xz}, E_{yz}\}$.  

\begin{center}
    \renewcommand{\arraystretch}{1.25}
    \begin{tabular}{|c|c c|c|}
        \hline
        $X_1$ &  \textbf{Nontrivial Orbits} && \textbf{Sublattice Rank}\\ \hline \hline
        $(xy)$ & $\{\alpha, \bar\alpha\}$, $\{x,y\}$, &$\{C_k(m,n), C_k(m,n+2)\}$  & 8 \\ \hline
        $(xzyw)$ & 
        $\{C_{xz}, C_{yz}, C_{xw}, C_{yw}\}$,& 
        $\{C_{\alpha z}, C_{\bar\alpha w}, C_{y \beta}, C_{x \bar\beta}\}$,& 6\\  
        &$\{C_{\bar\alpha z}, C_{\alpha w}, C_{y \bar\beta}, C_{x \beta}\}$,&
        $\{C_{\alpha \beta}, C_{\alpha \bar\beta}, C_{\bar\alpha \beta}, C_{\bar\alpha \bar\beta}\}$,&\\
        &$\{C_2(0,0), C_2(3,1)\}$,& 
        $\{C_2(0,1), C_2(3,0)\}$,&\\
        &$\{C_2(0,2), C_2(3,3)\}$,&
        $\{C_2(0,3), C_2(3,2)\},$&\\
        &$\{C_2(1,0), C_2(2,1)\}$,&
        $\{C_2(1,1), C_2(2,0)\}$,&\\
        &$\{C_2(1,2), C_2(2,3)\}$,&
        $\{C_2(1,3), C_2(2,2)\}$,&\\
        &$\{C_3(m+d, n+e)\}_{d,e\in\{0,2\}}$&   & \\ \hline
    \end{tabular}
\end{center}
\bigskip

\begin{center}
     \renewcommand{\arraystretch}{1.20}
    \begin{tabular}{|c|c c c|c|}
        \hline
        $X_4$ &  \textbf{Nontrivial Orbits} &&& \textbf{Sublattice Rank}\\ \hline \hline
        $(xy)$ & $\{\alpha, \bar\alpha\}$, $\{x,y\}$,& $\{C_k(0,0), C_k(0,2)\}$ && 12 \\ \hline
        $(xzyw)$ &  
        $\{C_{xz}, C_{yz}, C_{xw}, C_{yw}\}$,& 
        $\{C_{\alpha z}, C_{\bar\alpha w}, C_{y \beta}, C_{x \bar\beta}\}$,&  
        $\{C_{\bar\alpha z}, C_{\alpha w}, C_{y \bar\beta}, C_{x \beta}\}$,&6\\ 
        &$\{C_{\alpha \beta\sim \bar\alpha \bar\beta}, C_{\alpha \bar\beta\sim \bar\alpha \beta}\}$,&
        $\{C_2(0,0), C_2(0,2)\}$,& 
        $\{C_3(0,0), C_3(0,2)\}$,&\\ 
        &$\{E_x,E_y,E_z,E_w\}$, 
        &$\{E_{\alpha\sim\bar\alpha, E_\beta\sim\bar\beta}\}$, 
        &$\{E_{zx},E_{yz}, E_{wy}, E_{xw}\}$,&\\ 
        &$\{E_{zy},E_{yw}, E_{wx}, E_{xz}\}$&&& \\ \hline
    \end{tabular}
\end{center}
\bigskip

Surprisingly, we see that the invariant lattice of $(xy)$ on $X_1$, which we will denote by $S(xy)_1$, has rank 8, and the invariant lattice for $(xy)$ on $X_4$, which we denote by $S(xy)_4$, has rank 12. From Equation~\eqref{eq:mirrorRank}, we see these are possibly mirror lattices, so we need some further investigation. Unfortunately for $X_1$, the lattice generated by orbit sums is not the full invariant lattice. As described in Section~\ref{ss:latticebackground}, we need to describe the topology of the fixed locus of each automorphism. 

On $X_1$, we can see that the only points fixed by $(xy)$ are along the curve $D_1$ defined by $x-y=0$. We can view $D_1$ as the intersection of $V(x-y)$ in $\CC\PP^3$ with our K3 surface $X_1$. But in $\CC\PP^3$, the hypersurface $V(x-y)$ is isomorphic to $\CC\PP^2$, so we can think of $D_1$ as a degree 4 curve in $\CC\PP^2$. 
A degree $d$ curve in $\CC\PP^2$ has genus $\frac{(d-1)(d-2)}{2}$ (see e.g. \cite[Chapter V]{Hartshorne}). So in our case, we have a genus 3 curve fixed, and no other curves. Consulting Theorem~\ref{thm:nikulingark} and Figure~\ref{fig:2elementary}, we have $r=8$, $a=8$ and $\delta=1$.

On $X_4$ finding the fixed locus is slightly more complicated. Let us denote by $D_4$, the curve that is the image of $D_1$ under the automorphism $\sigma_L-\sigma_R$. It turns out the preimage of $D$ is two curves, $D_1$ and $V(x+y)\cap X_1$. Restricting $\sigma_L-\sigma_R$ to $D_1$, we obtain a degree 2 map $D_1\to D_4$ with four fixed points. We can use the Riemann--Hurwitz formula to compute the genus of $D_4$ as follows:
\[
\chi(D_1)=2\chi(D_4)-4. 
\]
As the genus of $D_1$ is $g=3$, we have $\chi(D_1)=2-2g_1=-4$, and so $\chi(D_4)=2-2g_4=0$, so the genus $g_4$ of $D_4$ is 1. 

Under the automorphism $(xy)$ on $X_4$, we have three fixed curves, namely $D_4$, $C_{z,\alpha\sim\bar\alpha}$, and $C_{w,\alpha\sim\bar\alpha}$. Hence, we have $g=1$, $k=2$. Consulting Theorem~\ref{thm:nikulingark} and Figure~\ref{fig:2elementary}, we have $r=12$, $a=8$, and $\delta=1$. 
Thus, we see that $S(xy)_1$ and $S(xy)_4$ are mirror lattices, and so $X_1$ and $X_4$ are LPK3 mirror surfaces.

\subsection{Self-Mirror Surface $X_2$}

\subsubsection{Geometry}

Similar to $X_4$, we use the lines on $X_1$ to derive a set of lines on $X_2$. Defining
\[
\sigma_2 := \left(\frac{1}{2},\frac{1}{2},0, 0\right)=4\sigma_L~,
\]
we get $\overline{G}_2=G_2/J_W = \langle \sigma_2 \rangle$ (and in the quotient $4\sigma_L=4\sigma_R$). Again, it suffices to consider the action of $\sigma_2$ on the lines on $X_1$. 

Under $\sigma_2$, we can see that all lines in family $C_1$ are invariant. 
Next, considering the action of $\sigma_2$ on the families $C_2$ and  $C_{3}$, we find the following sets of curves are identified:
\begin{align*}
    C_2(m,n) &\leftrightarrow C_2(m+2,n+2)\\
    C_3(m,n) &\leftrightarrow C_3(m+2,n+2).
\end{align*}

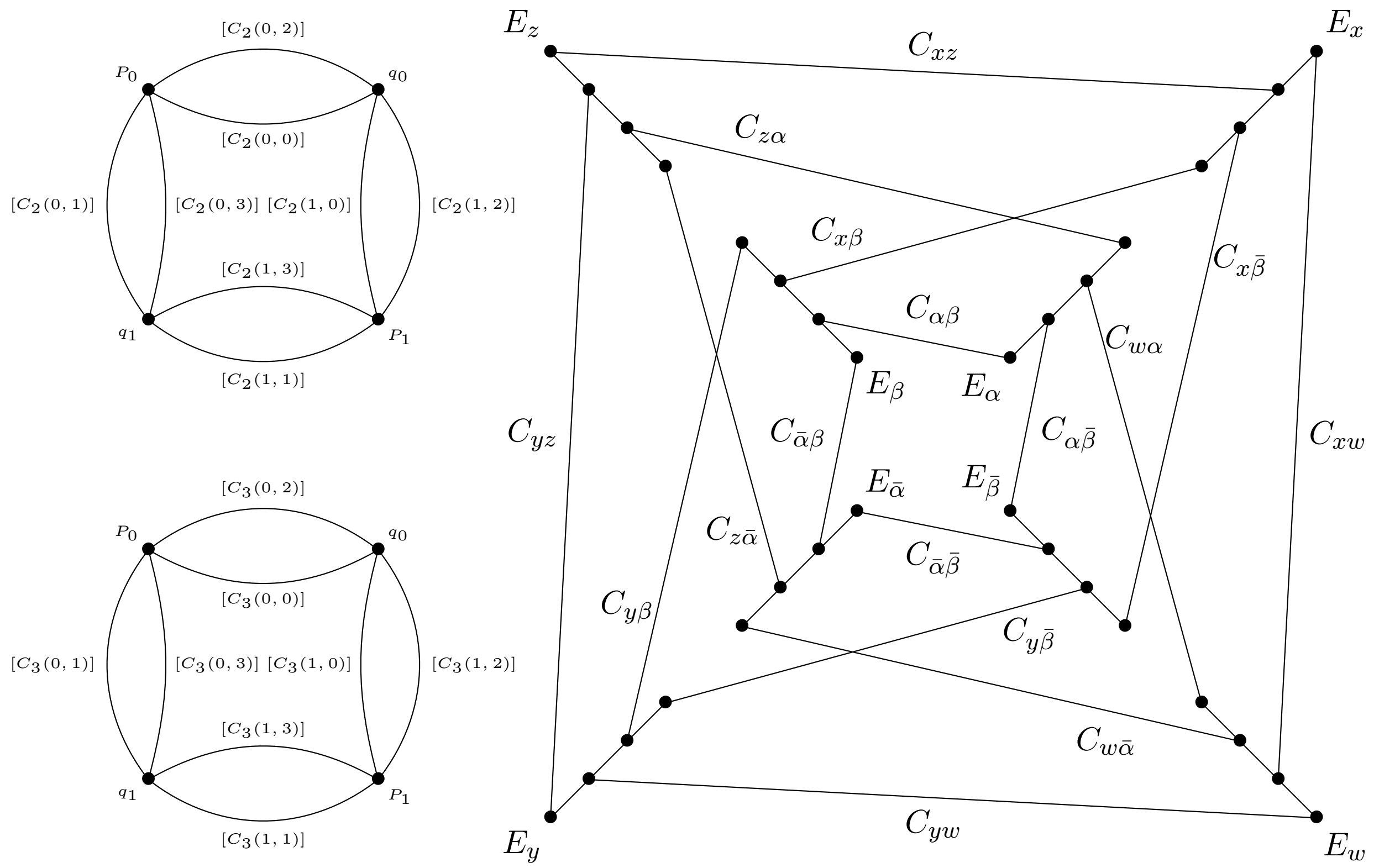
\begin{figure}[ht]
\centering\resizebox{\textwidth}{!}{\begin{tikzpicture}[xscale=.45,yscale=.45]
\tikzset{clbl/.style={font=\fontsize{4.5}{5.5}\selectfont}}

\begin{scope}[shift={(0,0)}]
\foreach \Point in {(-10,10),(-9,9),(-8,8),(-7,7),(-5,5),(-4,4),(-3,3),(-2, 2),(10,10),(9,9),(8,8),(7,7),(5,5),(4,4),(3,3),(2,2),(10,-10),(9,-9),(8,-8),(7,-7),(5,-5),(4,-4),(3,-3),(2,-2),(-10,-10),(-9,-9),(-8,-8),(-7,-7),(-5,-5),(-4,-4),(-3,-3),(-2,-2)}{
    \node at \Point {\textbullet};
}
\node at (-10.75,10.75){$E_z$};
\draw (-10,10)--(-7,7);
\node at (-1.25, 1.25){$E_\beta$};
\draw (-5,5)--(-2,2);
\node [above] at (0,9.5){$C_{xz}$};
\draw (-10,10)--(9,9);
\node at (-4.5, 8){$C_{z\alpha}$};
\draw (-8,8)--(5,5);
\node at (-2.5, 5.25){$C_{x\beta}$};
\draw (-4,4)--(7,7);
\node [above] at (0,2.5){$C_{\alpha\beta}$};
\draw (-3,3)--(2,2);

\node at (10.75,10.75){$E_x$};
\draw (10,10)--(7,7);
\node at (1.25, 1.25){$E_\alpha$};
\draw (5,5)--(2,2);
\node [right] at (9.5,0){$C_{xw}$};
\draw (10,10)--(9,-9);
\node at (8, 4.5){$C_{x\bar{\beta}}$};
\draw (8,8)--(5,-5);
\node at (5.25, 2.5){$C_{w\alpha}$};
\draw (4,4)--(7,-7);
\node [right] at (2.5,0){$C_{\alpha\bar{\beta}}$};
\draw (3,3)--(2,-2);

\node at (10.75,-10.75){$E_w$};
\draw (10,-10)--(7,-7);
\node at (1.25, -1.25){$E_{\bar{\beta}}$};
\draw (5,-5)--(2,-2);
\node [below] at (0,-9.5){$C_{yw}$};
\draw (10,-10)--(-9,-9);
\node at (4.5, -8){$C_{w\bar{\alpha}}$};
\draw (8,-8)--(-5,-5);
\node at (2.5, -5.25){$C_{y\bar{\beta}}$};
\draw (4,-4)--(-7,-7);
\node [below] at (0,-2.5){$C_{\bar{\alpha}\bar{\beta}}$};
\draw (3,-3)--(-2,-2);

\node at (-10.75,-10.75){$E_y$};
\draw (-10,-10)--(-7,-7);
\node at (-1.25, -1.25){$E_{\bar{\alpha}}$};
\draw (-5,-5)--(-2,-2);
\node [left] at (-9.5,0){$C_{yz}$};
\draw (-10,-10)--(-9,9);
\node at (-8,-4.5){$C_{y\beta}$};
\draw (-8,-8)--(-5,5);
\node at (-5.25, -2.5){$C_{z\bar{\alpha}}$};
\draw (-4,-4)--(-7,7);
\node [left] at (-2.5,0){$C_{\bar{\alpha}\beta}$};
\draw (-3,-3)--(-2,2);
\end{scope}

\coordinate (P0) at (-20.5,9);  \coordinate (q0) at (-14.5,9);
\coordinate (q1) at (-20.5,3);  \coordinate (P1) at (-14.5,3);
\draw (P0) to[bend left=38]  node[above,clbl]             {$[C_2(0,2)]$} (q0);
\draw (P0) to[bend right=30] node[below=-1pt,clbl]        {$[C_2(0,0)]$} (q0);
\draw (q1) to[bend left=30]  node[above=-1pt,clbl]        {$[C_2(1,3)]$} (P1);
\draw (q1) to[bend right=38] node[below,clbl]             {$[C_2(1,1)]$} (P1);
\draw (P0) to[bend right=38] node[left,clbl]              {$[C_2(0,1)]$} (q1);
\draw (P0) to[bend left=15]  node[pos=.5,right=-1pt,clbl] {$[C_2(0,3)]$} (q1);
\draw (q0) to[bend right=15] node[pos=.5,left=-1pt,clbl]  {$[C_2(1,0)]$} (P1);
\draw (q0) to[bend left=38]  node[right,clbl]             {$[C_2(1,2)]$} (P1);
\foreach \Point in {(P0),(q0),(q1),(P1)}{ \node at \Point {\textbullet}; }
\node [above left=-1pt,clbl]  at (P0) {$P_0$};
\node [above right=-1pt,clbl] at (q0) {$q_0$};
\node [below left=-1pt,clbl]  at (q1) {$q_1$};
\node [below right=-1pt,clbl] at (P1) {$P_1$};

\coordinate (P0b) at (-20.5,-3);  \coordinate (q0b) at (-14.5,-3);
\coordinate (q1b) at (-20.5,-9);  \coordinate (P1b) at (-14.5,-9);
\draw (P0b) to[bend left=38]  node[above,clbl]             {$[C_3(0,2)]$} (q0b);
\draw (P0b) to[bend right=30] node[below=-1pt,clbl]        {$[C_3(0,0)]$} (q0b);
\draw (q1b) to[bend left=30]  node[above=-1pt,clbl]        {$[C_3(1,3)]$} (P1b);
\draw (q1b) to[bend right=38] node[below,clbl]             {$[C_3(1,1)]$} (P1b);
\draw (P0b) to[bend right=38] node[left,clbl]              {$[C_3(0,1)]$} (q1b);
\draw (P0b) to[bend left=15]  node[pos=.5,right=-1pt,clbl] {$[C_3(0,3)]$} (q1b);
\draw (q0b) to[bend right=15] node[pos=.5,left=-1pt,clbl]  {$[C_3(1,0)]$} (P1b);
\draw (q0b) to[bend left=38]  node[right,clbl]             {$[C_3(1,2)]$} (P1b);
\foreach \Point in {(P0b),(q0b),(q1b),(P1b)}{ \node at \Point {\textbullet}; }
\node [above left=-1pt,clbl]  at (P0b) {$P_0$};
\node [above right=-1pt,clbl] at (q0b) {$q_0$};
\node [below left=-1pt,clbl]  at (q1b) {$q_1$};
\node [below right=-1pt,clbl] at (P1b) {$P_1$};

\end{tikzpicture}}
\caption{SLICE 1: lines on $X_2$. Intersections indicated with ``\textbullet''.}
\label{fig:X_2-lines-slice1}
\end{figure}

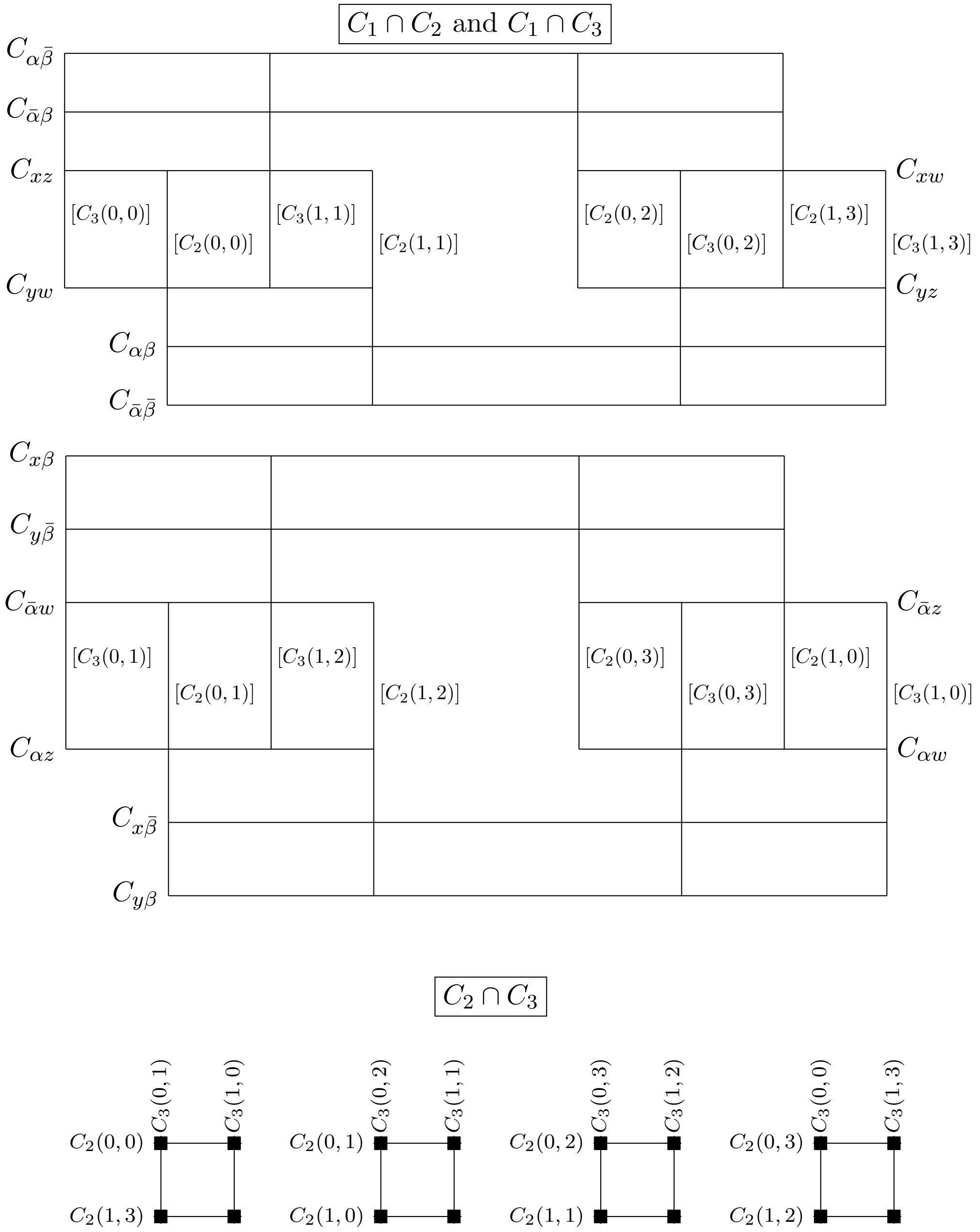
\begin{figure}[ht]
\centering

\begin{tikzpicture}[xscale=.7,yscale=.4, lab/.style={anchor=west, inner sep = 2pt}]

\node[draw,inner sep=3pt] at (0,7)  {$C_1\cap C_2$ and $C_1\cap C_3$};

\draw (-8,2)--(-2,2);
\node [left] at (-8,2){$C_{xz}$};
\draw (-8,-2)--(-2,-2);
\node [left] at (-8,-2){$C_{yw}$};
\draw (2,2)--(8,2);
\node [right] at (8,2){$C_{xw}$};
\draw (2, -2)--(8,-2);
\node [right] at (8,-2){$C_{yz}$};

\draw (-8,6)--(6,6);
\node [left] at (-8,6){$C_{\alpha\bar{\beta}}$};
\draw (-8,4)--(6,4);
\node [left] at (-8,4){$C_{\bar{\alpha}\beta}$};
\draw (-6,-6)--(8,-6);
\node [left] at (-6,-6){$C_{\bar{\alpha}\bar{\beta}}$};
\draw (-6,-4)--(8,-4);
\node [left] at (-6,-4){$C_{\alpha\beta}$};

\draw (-8,6)--(-8,-2);
\node [lab] at (-8, .5){\scalebox{0.8}{\small$[C_3(0,0)]$}};
\draw (-6,2)--(-6,-6);
\node [lab] at (-6, -.5){\scalebox{0.8}{\small$[C_2(0,0)]$}};
\draw (-4,6)--(-4,-2);
\node [lab] at (-4, .5){\scalebox{0.8}{\small$[C_3(1,1)]$}};
\draw (-2,2)--(-2,-6);
\node [lab] at (-2, -.5){\scalebox{0.8}{\small$[C_2(1,1)]$}};

\draw (2,6)--(2,-2);
\node [lab] at (2, .5){\scalebox{0.8}{\small$[C_2(0,2)]$}};
\draw (4,2)--(4,-6);
\node [lab] at (4, -.5){\scalebox{0.8}{\small$[C_3(0,2)]$}};
\draw (6,6)--(6,-2);
\node [lab] at (6, .5){\scalebox{0.8}{\small$[C_2(1,3)]$}};
\draw (8,2)--(8,-6);
\node [lab] at (8, -.5) {\scalebox{0.8}{\small$[C_3(1,3)]$}};

\end{tikzpicture}

\centering

\begin{tikzpicture}[xscale=.7,yscale=.5, lab/.style={anchor = west, inner sep = 2pt}]

\draw (-8,2)--(-2,2);
\node [left] at (-8,2){$C_{\bar\alpha w}$};
\draw (-8,-2)--(-2,-2);
\node [left] at (-8,-2){$C_{\alpha z}$};
\draw (2,2)--(8,2);
\node [right] at (8,2){$C_{\bar \alpha z}$};
\draw (2, -2)--(8,-2);
\node [right] at (8,-2){$C_{\alpha w}$};

\draw (-8,6)--(6,6);
\node [left] at (-8,6){$C_{x{\beta}}$};
\draw (-8,4)--(6,4);
\node [left] at (-8,4){$C_{y\bar{\beta}}$};
\draw (-6,-6)--(8,-6);
\node [left] at (-6,-6){$C_{y{\beta}}$};
\draw (-6,-4)--(8,-4);
\node [left] at (-6,-4){$C_{x\bar\beta}$};

\draw (-8,6)--(-8,-2);
\node [lab] at (-8, .5){\scalebox{0.8}{\small$[C_3(0,1)]$}};
\draw (-6,2)--(-6,-6);
\node [lab] at (-6, -.5){\scalebox{0.8}{\small$[C_2(0,1)]$}};
\draw (-4,6)--(-4,-2);
\node [lab] at (-4, .5){\scalebox{0.8}{\small$[C_3(1,2)]$}};
\draw (-2,2)--(-2,-6);
\node [lab] at (-2, -.5){\scalebox{0.8}{\small$[C_2(1,2)]$}};

\draw (2,6)--(2,-2);
\node [lab] at (2, .5){\scalebox{0.8}{\small$[C_2(0,3)]$}};
\draw (4,2)--(4,-6);
\node [lab] at (4, -.5){\scalebox{0.8}{\small$[C_3(0,3)]$}};
\draw (6,6)--(6,-2);
\node [lab] at (6, .5){\scalebox{0.8}{\small$[C_2(1,0)]$}};
\draw (8,2)--(8,-6);
\node [lab] at (8, -.5){\scalebox{0.8}{\small$[C_3(1,0)]$}};

\end{tikzpicture}


\centering

\newcommand{\twogrid}[4]{%

   \draw[step=2.0,black,thin] (-0.2,-0.2) grid (2.2,2.2);

  \foreach \a/\b [count=\k from 0] in {#3}{%
   
    \node[anchor=east,font=\scriptsize] at (-0.18,2-2*\k)
      {$#1(\a,\b)$};
  }

  \foreach \a/\b [count=\k from 0] in {#4}{%
      \node[
      rotate=90,
      font=\scriptsize
    ] at (2*\k,3.3) {$#2(\a,\b)$};
  }
  \foreach \Point in {(0,0), (0,2), (2,0), (2,2)}{
    \node at \Point {\tiny$\blacksquare$};
}
}

\vspace{.75cm}

\begin{tikzpicture}[x=.5cm,y=.5cm]

 \node[draw,inner sep=3pt] at (0,6) {$C_2\cap C_3$};

 \begin{scope}[shift={(-9,0)}]
 \twogrid{C_2}{C_3}
   {0/0,1/3}
   {0/1,1/0}
 \end{scope}

 \begin{scope}[shift={(-3,0)}]
 \twogrid{C_2}{C_3}
   {0/1,1/0}
   {0/2,1/1}
 \end{scope}

 \begin{scope}[shift={(3,0)}]
 \twogrid{C_2}{C_3}
   {0/2,1/1}
   {0/3,1/2}
 \end{scope}

 \begin{scope}[shift={(9,0)}]
 \twogrid{C_2}{C_3}
   {0/3,1/2}
   {0/0,1/3}
 \end{scope}

\end{tikzpicture}


\caption{SLICE 2 of lines on $X_2$. Intersections are depicted by intersections of the corresponding line in the diagram. Intersections depicted with {\tiny$\blacksquare$} indicate that the corresponding lines intersect in two (distinct) points.}
\label{fig:X_2-lines-slice2}
\end{figure}

Referring to Figure \ref{fig:X_1-lines-slice1}, there are eight $A_1$ singularities at the points labeled $t_x, t_y, t_z, t_w, t_\alpha, t_{\bar{\alpha}}, t_\beta, t_{\bar{\beta}}$. Since $\overline{G_2} \cong \ZZ/(2\ZZ)$, by \cite{Xiao} these are all the singularities resulting from this group action. Moreover, the families $C_2$ and $C_{3}$ do not intersect these points.
After resolving the singularities, we arrive at the diagram of curves in Figures \ref{fig:X_2-lines-slice1} and \ref{fig:X_2-lines-slice2}. 

As before, since the graph in Figure~\ref{fig:X_2-lines-slice1} is not planar, the intersections on $X_2$ are indicated by a ``\textbullet''. If two lines in the digram intersect without ``\textbullet'', there is no intersection between those two lines on $X_2$. In Figure~\ref{fig:X_2-lines-slice2}, however, any intersection between lines in the diagram \emph{does} indicate an intersection between the corresponding lines on $X_2$. Furthermore, the markings $\blacksquare$ indicate the the two lines intersecting there actually intersect in two distinct points on $X_2$. 

\subsubsection{Diagonal Automorphisms}

The diagonal automorphisms on $X_2$ are
\[
\overline{\Gdiag_2} := \Gdiag / G_2,\textrm{ where}\quad G_2 =\left\langle j_W, \sigma_2\right\rangle~.
\]
Again, we have $\overline{\Gdiag_2}$ generated by $\tau$ and $\sigma_R$, but now $4\sigma_R=4\sigma_L$ is trivial, so $\sigma_R$ has order 4. 
We find
\[\overline{\Gdiag_2} \cong \ZZ_2 \oplus \ZZ_4\]
given by 
\[
(m\tau + n\sigma_R) + G_2 \mapsto (m, n)~.
\]

Adjusting the previous table to this new group, we again show which diagonal symmetries are nonsymplectic: 

\begin{center}
    \renewcommand{\arraystretch}{1.25}
    \begin{tabular}{|c|c|c||c|c|c|}
        \hline
        $\bm{(m, n)}$ & Order & Determinant & $\bm{(m, n)}$ & Order & Determinant\\ \hline \hline
        (0, 0) & 1 & 1 & (1,0) & 2 & $-1$ \\ \hline
        (0, 1) & 8 & $i$ & (1,1) & 8 & $-i$ \\ \hline
        (0, 2) & 4 & $-1$ & (1,2) & 4 & $1$ \\ \hline
        (0, 3) & 8 & $-i$ & (1,3) & 8 & $i$ \\ \hline
    \end{tabular}
\end{center}

\subsubsection{Diagonal Invariant Sublattices on $X_2$}

We proceed in a similar way as before. From \eqref{eq:sigmaR} and \eqref{eq:sigmaL}, we can quickly compute the orbits of each nonsymplectic diagonal automorphism. 

For nonsymplectic automorphisms, we must only consider the following: 
\begin{itemize}
    \item First $(0,1)$ and $(1,1)$ are purely nonsymplectic of order 4; 
    \item Second, $(0,2)$ and $(1,0)$ are nonsymplectic involutions. 
\end{itemize}

As before $\{\beta,\bar\beta\}$ represents all of the orbits with $\{C_{\ell\beta},C_{\ell,\bar\beta}\}$, but also $\{E_\beta,E_{\bar\beta}\}$. 

\begin{center}
    \renewcommand{\arraystretch}{1.25}
    \begin{tabular}{|c|l|c|}
        \hline
        $\bm{(m, n)}$ &  \textbf{Nontrivial Orbits} & \textbf{Sublattice Rank}\\ \hline \hline
        (0, 1) &  $\{\beta, \bar{\beta}\}, \{C_2(0+i,d+i), C_3(0+i,d+i)\}_{i\in\{0,1\}, d\in \ZZ/4\ZZ}$ & 14 \\ \hline
        (0, 2) &  $\{C_k(0+i,d+i)\}_{i\in\{0,1\}, d\in \ZZ/4\ZZ}$ & 18 \\ \hline
        (1, 0) &  $\{\alpha, \bar{\alpha}\}, \{\beta, \bar{\beta}\}, \{C_2(0,d), C_2(d,0)\}_{d\neq 0}$ & 14 \\ \hline
        (1, 1) &  $\{\alpha, \bar{\alpha}\}, \{C_2(0+i,d+i), C_3(d+i,0+i)\}_{i\in\{0,1\}, d\in \ZZ/4\ZZ}$ & 14 \\ \hline
    \end{tabular}
\end{center}
\bigskip
Thus we can see that no diagonal symmetry will suffice to show that this surface is LPK3 mirror to itself.

\subsubsection{Non-Diagonal Automorphisms}

From Section~\ref{ss:nondiagonal}, we also have the same nondiagonal symmetries on $X_2$ as we did on the other two surfaces induced by permuting the variables. As before, four of these are nonsymplectic, but again we only consider $(xzyw)$ and $(xy)$. 

We can compute the ranks of the invariant lattice via the orbit sums. We represent all of the orbits of the form $\{C_{\alpha \ell}, C_{\bar\alpha \ell}\}$, and similarly with $\{x,y\}$ as before. This includes on $X_4$ orbits containing exceptional curves, such as $\{E_{x}, E_{y}\}$ and $\{E_{\alpha}, E_{\bar\alpha}\}$. These are given in the following table: 

\begin{center}
    \renewcommand{\arraystretch}{1.25}
    \begin{tabular}{|c|c c|c|}
        \hline
        $X_2$ &  \textbf{Nontrivial Orbits}& & \textbf{Sublattice Rank}\\ \hline \hline
        $(xy)$ & $\{\alpha, \bar\alpha\}$, $\{x,y\}$, &$\{C_k(m,n), C_k(m,n+2)\}$ & 10\\ \hline
        $(xzyw)$ &  
    
        $\{C_{xz}, C_{yz}, C_{xw}, C_{yw}\}$,& 
        $\{C_{\alpha z}, C_{\bar\alpha w}, C_{y \beta}, C_{x \bar\beta}\}$,&6\\  
        &$\{C_{\bar\alpha z}, C_{\alpha w}, C_{y \bar\beta}, C_{x \beta}\}$,& 
        $\{C_{\alpha \beta}, C_{\alpha \bar\beta}, C_{\bar\alpha \beta}, C_{\bar\alpha \bar\beta}\}$,&\\
        &$\{E_x,E_y,E_z,E_w\}$,& 
        $\{E_\alpha, E_{\bar\alpha}, E_\beta, E_{\bar\beta}\}$&\\
        &$\{C_2(0,0), C_2(1,3)\}$,& 
        $\{C_2(0,2), C_2(1,1)\}$,&\\
        &$\{C_2(1,0), C_2(0,3)\}$,&
        $\{C_2(1,2), C_2(0,1)\}$,&\\
        &$\{C_3(0,0), C_3(0,2)\}$,& 
        $\{C_3(0,1), C_3(1,0)\}$,&\\
        &$\{C_3(0,3), C_3(1,2)\}$,&
        $\{C_3(1,1), C_3(1,3)\}$&\\\hline
        
    \end{tabular}
\end{center}
\bigskip

Surprisingly, we see that the invariant lattice of $(xy)$ has rank 10, and so by Equation~\eqref{eq:mirrorRank}, this could be its own mirror lattice. Recall that the mirror lattice of a 2-elementary lattice is the 2-elementary lattice with invariants $20-r, a, \delta$, so any 2-elementary lattice is isomorphic to its mirror. So, there is nothing more to check, and we see that $X_2$ is LPK3 mirror to itself. It turns out the invariants are $g=1$, $k=0$, so we have $r=10$, $a=10$. We did not compute $\delta$.

\appendix 

\section{Transformation to Fermat}\label{appendix:A}

There exists a linear transformation on $\CC^4$ which takes $W$ from Equation \eqref{eq:definition-W} to a Fermat-type polynomial. Define
\begin{align*}
    c &:= \frac{1}{2^{1/4}} \\
    \zeta &:= \exp\left(\frac{\pi i}{4}\right)
\end{align*}
and linear map $B:\CC^4 \to \CC^4$ by
\[
B(x, y, z, w) :=
\begin{pmatrix}
    c & -c\zeta & 0 & 0 \\
    c & c\zeta & 0 & 0 \\
    0 & 0 & 1 & 0 \\
    0 & 0 & 0 & 1
\end{pmatrix}
\begin{pmatrix}
    x \\
    y \\
    z \\
    w \\
\end{pmatrix}~.
\]

We can transform $W$ via $W'(\vec{x})=W(B\vec{x})$, whereupon $W'=x^4+y^4+z^3w+zw^3$. As $B$ is invertible, we get an isomorphism
\[
B^{-1}:V(W)\to V(W').
\]
Note we could do a similar thing with the last two variables, giving us an isomporphism between $X_1$ and the Fermat quartic. 

However, we would like to point out, that although BHK mirror symmetry has been studied for the Fermat quartic, we cannot leverage the BHK mirror of the Fermat together with this isomorphism to get a statement about the about BHK mirror symmetry for $X_1$, at least not in the naive way. 

Observe if $g\in G^{diag}_{W}$, then we obtain a symmetry of $W'$ via $g'=A^{-1}g A$. However, one can check that in general $g'$ is not a diagonal symmetry. BHK mirror symmetry depends heavily on the diagonal symmetry groups. 


This example demonstrates that to prove the equivalence of BHK and LPK3 in other cases, it is likely that non-diagonal automorphisms are required. Or, BHK needs to be extended to include K3 surfaces constructed from non-diagonal group actions.
 
 

\newpage

\bibliographystyle{amsplain}
\bibliography{references}

\end{document}